\documentclass[reqno,twoside,11pt,english]{amsart}
\usepackage{amsmath,amsfonts,amssymb,amsthm,epsfig}
\newcommand{\RNum}[1]{\uppercase\expandafter{\romannumeral #1\relax}}
\usepackage{comment}

\usepackage{color}
\usepackage[final,colorlinks,linkcolor=blue,anchorcolor=red,citecolor=blue]{hyperref}
\usepackage{graphicx}
\usepackage{titletoc,epsf}
\usepackage{amsmath,amsfonts,latexsym,amsthm,amsxtra,amssymb,bbm}
\allowdisplaybreaks
\usepackage{graphicx,float}

\allowdisplaybreaks

\let\lam=\lambda

\let\f=\frac

\let\Om=\Omega

\let\th=T
\let\pa=\partial

\def\R{\mathbf R}

\DeclareMathOperator{\supp}{supp}

\newtheorem{theorem}{Theorem}[section]
\newtheorem{lemma}[theorem]{Lemma}
\newtheorem{proposition}[theorem]{Proposition}
\newtheorem{corollary}[theorem]{Corollary}
\newtheorem{definition}[theorem]{Definition}
\newtheorem{remark}[theorem]{Remark}

\newtheorem{thmA}{Theorem}[section]

\newtheorem{Ques}[thmA]{Question}

\numberwithin{equation}{section}

\numberwithin{equation}{section}

\topmargin - .23 in

\begin{document}
	
	\title{Global weak solution of 3-D focusing energy-critical nonlinear Schr\"odinger equation}
	
	\author{Xing Cheng, Chang-Yu Guo and Yunrui Zheng}

	\address[X. Cheng]{School of mathematics, Hohai University, Nanjing 210098, Jiangsu, P. R. China}
	\email{{\tt chengx@hhu.edu.cn}}
	
	\address[C.-Y. Guo]{Research Center for Mathematics and Interdisciplinary Sciences, Shandong University, 266237 Qingdao, P. R. China, Frontiers Science Center for Nonlinear Expectations, Ministry of Education, P. R. China and Department of Physics and Mathematics, University of Eastern Finland
		80101 Joensuu, Finland}
	\email{{\tt changyu.guo@sdu.edu.cn}}
	
	\address[Y. Zheng]{ School of Mathematics, Shandong University, Shandong 250100, Jinan, P. R. China.}
	\email{{\tt yunrui\_zheng@sdu.edu.cn}}
	
	\dedicatory{Dedicated to Professor~Gongbao Li on the occasion of his 70th birthday}

		\subjclass[2020]{Primary: 35L05; Secondary: 35L71, 35Q40, 35P25, 35B40}
	\keywords{Nonlinear Schr\"odinger equation, Ginzburg-Landau equation, well-posedness, blow-up, weak solutions}
	

	\begin{abstract}
		In this article, we prove the existence of global weak solutions to the three-dimensional focusing energy-critical nonlinear Schr\"odinger (NLS) equation in the non-radial case, which solves an open problem of T. Cazenave \cite{C}. Furthermore, we prove the weak-strong uniqueness for some class of initial data. The main ingredient of our new approach is to use solutions of an energy-critical Ginzburg-Landau equation as approximations for the corresponding nonlinear Sch\"ordinger equation.
		
		In our proofs, we first show the dichotomy of global well-posedness versus finite time blow-up of energy-critical Ginzburg-Landau equation in $\dot{H}^1( \mathbb{R}^d)$ for $d = 3,4 $ when the energy is less than that of the stationary solution $W$.
		We follow the strategy of C. E. Kenig and F. Merle \cite{KM1,KM}, using a concentration-compactness/rigidity theorem argument to reduce the global well-posedness to the exclusion of a critical element.
		The critical element is ruled out by dissipation of the Ginzburg-Landau equation, including local smoothness, backwards uniqueness and unique continuation.
		The existence of global weak solution to the three dimensional focusing energy-critical nonlinear Schr\"odinger equation in the non-radial case then follows from the global well-posedness of the energy-critical Ginzburg-Landau equation via a limitation argument. We also adapt the arguments of M. Struwe \cite{Stw01,Stw02} to prove the weak-strong uniqueness when the $\dot{H}^1$-norm of the initial data is bounded by a constant depending on the stationary solution $W$.
		
%
%
	\end{abstract}
	\maketitle
	\tableofcontents
	
	\section{Introduction}\label{se1}
	
	The energy-critical nonlinear Schr\"odinger (NLS) equation
	\begin{align}\label{eq1.1v18}
		\begin{cases}
			i \partial_t v + \Delta v +   \mu  f(v) = 0 , \\
			v(0) = v_0 \in \dot{H}^1 ( \mathbb{R}^d),
		\end{cases}
	\end{align}
	where $d \ge 3$ is the dimension, $\mu = \pm 1$, $f(v) = |v|^\frac4{d-2 } v$, and $v: \mathbb{R}_t \times \mathbb{R}^d_x \to \mathbb{C}$, has been studied intensively in the last three decades.  If $\mu = 1$, it is focusing; if $\mu = -1$, it is defocusing. Its energy as Hamiltonian is
	\begin{align*}
		E(v) =  \int_{\mathbb{R}^d  } \frac12 |\nabla v|^2 -  \mu \frac{d- 2 }{2d }  |v|^\frac{2d}{d-2} \,\mathrm{d}x.
	\end{align*}
	Regarding the theory of energy-critical nonlinear Schr\"odinger equation \eqref{eq1.1v18}, there is a distinguished problem that has significantly promoted the development of the field in the last two decades, for which we formulate as follows.

	\begin{Ques}[Global well-posedness and scattering]\label{ques:A}
		\textit{Do strong solutions to \eqref{eq1.1v18} globally exist and necessarily scatter in $\dot{H}^1$? }
	\end{Ques}
	
	Recall that a function $v \in C_t^0 \dot{H}^1_x( I \times \mathbb{R}^d)$ is a strong-$\dot{H}^1$ solution of \eqref{eq1.1v18}, where $ I\subset \mathbb{R}$ is an interval containing $0$,  if it satisfies the Duhamel formula
	\begin{align*}
		v(t) = e^{i t \Delta} v_0  +   i \mu \int_0^t e^{i (t - s) \Delta}  f(v(s)) \,\mathrm{d}s,
	\end{align*}
	for all $t \in I$. We refer to the interval $I$ as the lifespan of $v$, and call $v$ a maximal-lifespan solution if $I$ cannot be extended to any strictly larger interval. The solution $v$ is global if $I = \mathbb{R}$.
	\medskip
	
	In the \textbf{defocusing case}, J. Bourgain \cite{B} first proved the global well-posedness and scattering for \emph{the radial initial data} when $d = 3, 4$, where he introduced the induction on energy method and spatial localized Morawetz estimate. M. Grillakis \cite{Gr} provided another proof for the global well-posedness. Later, T. Tao \cite{T} considered the global well-posedness and scattering for \emph{the radial initial data} in higher dimensions.
	
	For \emph{the non-radial initial data}, J. Colliander, M. Keel, G. Staffilani, H. Takaoka, and T. Tao \cite{CKSTT} made a major breakthrough and proved the global well-posedness and scattering when $d = 3$. Later, Ryckman and Visan \cite{RV} proved global well-posedness and scattering when $d = 4$, M. Visan \cite{V2} proved global well-posedness and scattering for $d \ge 5$. We also refer to \cite{CGZ,KV1,V1} for a shorter treatment of the global well-posedness and scattering for $d = 3,4,5$ using the long time Strichartz estimate.
	
	\medskip
	
	In the \textbf{focusing case}, C. Kenig and F. Merle \cite{KM} invented the concentration-compactness/rigidity theorem method to establish the global well-posedness and scattering when $d = 3, 4, 5$ for \emph{radial initial data}. Later, R. Killip and M. Visan \cite{KV2010} proved global well-posedness and scattering for \emph{the non-radial initial data} when $d \ge 5$. Recently, B. Dodson \cite{D2019} proved the global well-posedness and scattering for \emph{the non-radial initial data} when $d = 4$.
	\medskip
	
	We can summarize the above results in the following theorem:
	\begin{thmA}[GWP \& scattering of the energy-critical nonlinear Schr\"odinger equation] \label{con1.1v18}
		Let $v_0 \in \dot{H}^1$, when $\mu = - 1$, the corresponding solution $v$ to \eqref{eq1.1v18}
		is global for $d \ge 3$. When $\mu =  1$ and $E(v_0) < E(W)$, where
		\begin{align}\label{eq1.3v19}
			W(x) :=  \frac1{ \left( 1 + \frac{  |x|^2}{ d(d-2) }  \right)^\frac{d-2}2 },
		\end{align}
		the following conclusions hold:
		\begin{itemize}
			\item[(1)] If $\|\nabla v_0 \|_{L^2} < \| \nabla W\|_{L^2}$, then the solution $v$ of \eqref{eq1.1v18} is global for $d \ge 4$ and this is also true for $d = 3$ in the radial case,
			\item[(2)] If $\|\nabla v_0 \|_{L^2} >  \|\nabla W\|_{L^2}$, and if {either} $v_0 \in L^2\left(\mathbb{R}^d, |x|^2 \mathrm{d}x \right)$ or $v_0 \in H^1$ is radial, then the solution $v$ of \eqref{eq1.1v18} blows up in finite time for $d \ge 3$.
		\end{itemize}
		Furthermore, the global solution scatters in $\dot{H}^1$.
	\end{thmA}
	
	Thus, according to Theorem \ref{con1.1v18}, the global well-posedness 
	for 3-d non-radial focusing NLS \textbf{remains as an open question}. It seems to be very difficult and challenging  to prove the global well-posedness and scattering of the focusing energy-critical nonlinear Schr\"odinger equation \emph{in the non-radial case}. The most natural approach to solve this hard problem would be using the concentration-compactness/rigidity theorem argument invented by C. E. Kenig and F. Merle \cite{KM,KM1}, which has been proved to be a powerful tool when dealing with the non-radial problem in the subsequent important development in \cite{D2019,KV2010,KV1,V1} for $d\ge 4$. However, it seems difficult to adapt this approach to solve the problem when $d = 3$. In general, the minimal energy blowup solution behaves like the stationary solution $W(x)$ in \eqref{eq1.3v19}, which does not belong to $L^2(\mathbb{R}^d)$ when $d = 3,4$; see \cite{KV2010}. In \cite{D2019}, B. Dodson successfully proved that the minimal energy blowup solution belongs to certain weak $L^2$ space when $d=4$ which plays a vital role to prove the rigidity theorem when $d =4$. But it seems that the argument in \cite{D2019} cannot be  extended to the three-dimensional case because the minimal energy blowup solution did not belong to any $L^2$ type space. Therefore, new ideas are needed to handle the case $d = 3$.
	
	Somewhat surprisingly, even the existence of global weak solutions to \eqref{eq1.1v18} when $\mu = 1$ in three dimension has still been open, which was already proposed twenty years ago in \cite[Remark 9.4.8]{C}.
	\begin{Ques}[\cite{C}]\label{ques:C}
		\textit{Do weak solutions to \eqref{eq1.1v18} globally exist when $\mu = 1$ and $d=3$}?
	\end{Ques}
	
	The main motivation of this paper is to give a satisfied answer to Question \ref{ques:C}. For this, recall first the definition of weak solutions.
	\begin{definition}\label{def:weak solution}
	A function $v \in L_t^\infty {H}^1_x( I \times \mathbb{R}^d)$ is a weak-${H}^1$ solution of \eqref{eq1.1v18} on an interval  $I$ containing $0$ (called the lifespan of $v$) if it satisfies the Duhamel formula
	\begin{align*}
		v(t) = e^{i t \Delta} v_0 +  i \mu  \int_0^t e^{i (t - s) \Delta}  f(v(s)) \,\mathrm{d}s,
	\end{align*}
	for almost every $t \in I$. We say that $v$ is a maximal-lifespan solution if the lifespan $I$ of $v$ cannot be extended to any strictly larger interval, and $v$ is a global solution if $I = \mathbb{R}$.
	\end{definition}
	
	We would like to point that in many parabolic PDEs or PDEs from Fluid, global weak solutions only require certain equation holds on each compact interval $I\subset \R$, which is weaker than the definition above. Our definition of global weak solution is similar to the one introduced by T. Tao \cite{Tao-2009-APDE}, where he considered global existence of weak solutions to the focusing mass-critical NLS.
	
	The main result of this paper is to prove the global existence and weak-strong uniqueness of weak solutions to \eqref{eq1.1v18} in the non-radial case when $d=3$.
	\begin{theorem}[Global weak solution of the 3-dimensional focusing energy-critical NLS]\label{th1.2v21}
		Fix $v_0 \in H^1$ and let $W$ be the stationary solution of \eqref{eq1.1v18} when $\mu=1$ given by \eqref{eq1.3v19}. Then we have the following result.
		
		\begin{itemize}
			\item[1)] Global existence: {If} $E(v_0) < E(W)$ and $\|\nabla v_0 \|_{L^2} < \| \nabla W\|_{L^2}$, {then} there exists a global weak solution $v\in L_t^\infty H_x^1(\mathbb{R} \times \mathbb{R}^3 ) \bigcap C_t^0 L^2_x(\mathbb{R} \times \mathbb{R}^3 )$ of \eqref{eq1.1v18} satisfying the energy inequality
			\begin{align}\label{eq1.4v26}
				E(v(t)) \le E(v_0)
			\end{align}
			and mass conservation
			\begin{align*}
				M(v(t)) = M(v_0).
			\end{align*}
			
			\item[2)] Weak-strong uniqueness: There exists a universal constant $C= C \left(W\right)>0$ such that if $\tilde{v}  \in C_tH^2_x( \mathbb{R} \times \mathbb{R}^3 ) $ is the global strong solution to  \eqref{eq1.1v18} with the initial data $v_0 \in H^2(\mathbb{R}^3)$ satisfying $E(v_0) < E(W)$ and $\|\nabla v_0\|_{L^2} < \min \left\{C,\|\nabla W\|_{L^2} \right\}$ and if ${v}$ is the global weak solution to \eqref{eq1.1v18} with initial datum $v_0$, then $v\equiv \tilde{v} $.
			
		\end{itemize}
		
	\end{theorem}
	
	\begin{remark}
		{\upshape	1)} It would not be an easy task if one tries to prove Theorem \ref{th1.2v21} part 1) via the standard Galerkin method due to the loss of compactness in the entire space $\mathbb{R}^3$; see \cite[Remark 9.4.8]{C} and Section \ref{sec:concluding remark} below for more detailed explanation.
		
		{\upshape	2)} We expect that the construction of global weak solutions would be a valuable first try to solve the global existence for strong solutions. In particular, if one can relax the assumption  $\|\nabla v_0\|_{L^2} < \min \left\{C,\|\nabla W\|_{L^2} \right\}$ to the weaker one $\|\nabla v_0\|_{L^2} < \|\nabla W\|_{L^2}$ to obtain the weak-strong uniqueness as in the theorem, then this would imply that the global weak solution provided by Theorem \ref{th1.2v21} is indeed the only possible strong solution of \eqref{eq1.1v18}. Thus one could consider further regularity of the weak solution provided by Theorem \ref{th1.2v21} so as to obtain a  global strong solution.
	\end{remark}
	
	We refer the interested readers to \cite{C} for a proof of the existence of global weak solutions to the \emph{defocusing} nonlinear Schr\"odinger equation. It is, however, impossible to extend the proof there to the focusing case: in \cite{C}, he truncated the nonlinear term, which destroys the variational structure of the equation in the focusing case. To prove Theorem \ref{th1.2v21}, we shall turn to study the global well-posedness of the following  Ginzburg-Landau equations
	\begin{align*}
		\partial_t u = \gamma u + \left(a + i \alpha \right) \Delta u - \left(b + i \beta \right) |u|^{p-1} u,
	\end{align*}
	where $u(t,x)$ is complex-valued functions of $(t,x) \in \mathbb{R}_+\times \mathbb{R}^d $, $a, b, \alpha, \beta, \gamma \in \mathbb{R}$, and $p > 1$. This class of equations has drawn the attentions of many mathematicians and physicists, because of its close relation to many physical phenomenons in fluids and superconductivity, etc. Especially, the Ginzburg-Landau equations could be derived from the B\'enard convection and turbulence. We refer the interested readers to \cite{CH93,GL,GJL,Gr1} for a complete overview.
	
	It is well-known that the Ginzburg-Landau equation is closely related to the Schr\"odinger equation. For instance, the nonlinear Schr\"odinger equation is the inviscid limit of Ginzburg-Landau equation; see \cite{BJ,W0,HW,Cheng-Guo-Zheng-2024} and the references therein. Although there already exists many results on the well-posedness of the Ginzburg-Landau equation, to the best of our knowledge, the global well-posedness of this very important model when $b<0$ has not yet been studied in the literature.
	
	In this article, we shall consider the Cauchy problem of  the following energy-critical Ginzburg-Landau equation in $\mathbb{R}^d$ when $d \ge 3$,
	\begin{align}\label{eq1.1}
		\begin{cases}
			u_t -  z \Delta u -    z f(u)= 0, \\
			u(0,x) =  u_0(x),
		\end{cases}
	\end{align}
	where the complex number $z$ satisfies $|z|$=1 with $0 < \Re z \le 1$, and $f(u) =  |u|^\frac4{d-2}u$.
	Equation \eqref{eq1.1} is the $L^2$ gradient flow of the energy
	\begin{align*}
		E(u) = \int_{\mathbb{R}^d  } \frac12 |\nabla u|^2 -   \frac{d- 2 }{2d }  |u|^\frac{2d}{d-2} \,\mathrm{d}x.
	\end{align*}
	We now prove the dichotomy of global well-posedness versus finite time blow-up of \eqref{eq1.1} when $d=3, 4$, which is our second main result.
	\begin{theorem} \label{th1.1v3}
		Let $u_0 \in \dot{H}^1(\mathbb{R}^d )$ and $E(u_0) <  E(W)$, where $W$ is the stationary solution of \eqref{eq1.1} given by \eqref{eq1.3v19}.
		Then we have
		\begin{itemize}
			\item[(1)] When $ \|\nabla u_0 \|_{L^2} < \|\nabla W\|_{L^2}$, the solution $u$ of \eqref{eq1.1} is global and satisfies
			$\lim\limits_{t \to \infty} \|u(t) \|_{\dot{H}^1} = 0$.
			\item[(2)] When $\|\nabla u_0 \|_{L^2} > \|\nabla W\|_{L^2}$ and $u_0 \in H^1(\mathbb{R}^d )$, the solution $u$ of \eqref{eq1.1} blows up forward in finite time.
			\item[(3)] The case $\|\nabla u_0 \|_{L^2} = \|\nabla W\|_{L^2}$ cannot happen.
		\end{itemize}
		
	\end{theorem}
	
	The stationary solution $W$ in \eqref{eq1.3v19} provides a threshold for global result and blow-up criterion,
	which plays the same role as ground state of \eqref{eq1.1v18} in the focusing case. The tools we use to prove Theorem \ref{th1.1v3} are those dealing with the energy-critical nonlinear dispersive equations \cite{KM,KM1,D2019,KV2010} as well as parabolic equations \cite{GR,KK11}, especially using the method of concentration-compactness/rigidity theorem developed in \cite{KM,KM1}. Note that when $\Re z $ is close to zero, the Ginzburg-Landau equation has similar behavior to the nonlinear Schr\"odinger equation; when $\Im z $ is close to zero, the Ginzburg-Landau equation has similar behavior to the heat equation. So our strategy to prove Theorem \ref{th1.1v3} is that we first use the concentration-compactness argument to construct the critical element or minimal energy blowup solution, and then we use the rigidity established on the effect of dissipation to exclude the possibility of blowup in finite time. To prove the existence of the critical element, we derive the linear profile decomposition of the Ginzburg-Landau equation by following the argument in \cite{G,GKP,Ke,BG}, and then use the method in \cite{KM,KV2010,GR,KK11} to establish the existence of the critical element. This procedure holds for all $d\ge 3$.
	To prove the rigidity, we employ some parabolic tools, such as backwards uniqueness and unique continuation developed in \cite{ESS032} for Navier-Stokes equations and local smallness regularity criterion. However, the rigidity only holds for $d=3, 4$.
	
	With the aid of Theorem \ref{th1.1v3}, we are able to prove, using a direct compactness argument, the existence of global weak solutions to the three dimensional focusing energy-critical nonlinear Schr\"odinger equation in the non-radial case. Motivated by the argument in \cite{Stw01,Stw02}, where M. Struwe proved the uniqueness for the weak solution of the defocusing energy critical nonlinear wave and Schr\"odinger equations, we  provide a weak-strong uniqueness result of the three dimensional focusing energy-critical nonlinear Schr\"odinger equation when the $\dot{H}^1$-norm of the initial data is bounded by a constant depending on $W$.
	
	We would like to point out that the use of solutions of (real or complex) Ginzburg-Landau type equation as proper approximations of geometric models was widely used in geometric analysis. For instance, in the seminal work of Y. M. Chen and M. Struwe \cite{Chen-Struwe-1989}, the authors used solutions of (real) Ginzburg-Landau type equations to approximate solutions of the harmonic mapping heat flow and established global existence of weak solutions to the heat flow. This approximation argument was later extended by Y. M. Chen and F.-H. Lin \cite{Chen-Lin-1993} to study evolution of harmonic maps with Dirichlet boundary conditions, and by Struwe \cite{Struwe-1991} to consider the evolution of harmonic mappings with free boundaries. Similar ideas were used by C. Y. Wang \cite{Wang-2002} (in the time-independent case) to study convergence properties of $p$-harmonic maps. In another seminal work \cite{Lin-Riviere-1999}, F.-H. Lin and T. Rivi\`ere used the (complex) Ginzburg-Landau relation procedure as approximation of unimodular harmonic maps.
	
	The structure of the paper is as follows. In Section \ref{se2v23}, we present the local well-posedness theory and also some variational  facts of \eqref{eq1.1}. In Section \ref{se3v23}, we show the finite time blowup part of Theorem \ref{th1.1v3}.  In Section \ref{se5v23}, we give the existence of the critical element. In Section \ref{se6v23}, we prove the rigidity theorem. In Section \ref{se4v23}, we show the asymptotic decay of global solutions. In Section \ref{se7v23}, we present the proof of Theorem \ref{th1.2v21}. In the final section, Section \ref{sec:concluding remark}, we shall point out that the standard Galerkin method cannot be easily adapted to prove Theorem \ref{th1.2v21}.
	
	Our notations are standard. We use $X \lesssim Y$ when $X \le CY$ for some constant $C> 0$ and $X \sim Y$ when $X \lesssim Y \lesssim X$.
	
	For any space-time slab $I \times \mathbb{R}^d$, we use $L_t^q L_x^r(I \times \mathbb{R}^d)$ to denote the space of functions $u: I \times \mathbb{R}^d \to \mathbb{C}$, whose norm is
	\begin{align*}
		\|u\|_{L_t^q L_x^r (I \times \mathbb{R}^d)} : = \left( \int_I \|u(t) \|_{L_x^r(\mathbb{R}^d) }^q \,\mathrm{d}t \right)^\frac1q < \infty,
	\end{align*}
	with the usual modification when $q$ or $r$ are equal to infinity. When $q = r$, we abbreviate $L_t^q L_x^q$ as $L_{t,x}^q$.
	
	We define the $S(I)$-norm of $u$ on the time slab $I$ to be
	\begin{align*}
		\|u\|_{S(I)} : = \left( \int_I \int_{\mathbb{R}^d } |u(t,x)|^\frac{2(d+2)}{d-2}  \,\mathrm{d}x \mathrm{d}t \right)^\frac{d-2}{2(d+2)} .
	\end{align*}
	
	\section{Auxiliary results} \label{se2v23}

	In this section, {we recall the local well-posedness theory of \eqref{eq1.1}, whose proofs can be found} in \cite{GV2,GJL,HW,W0}. We also establish some variational results of \eqref{eq1.1}.

	\subsection{Local well-posedness theory}
	
	We first give the definition of strong $\dot{H}^1$ solution of \eqref{eq1.1}.
	\begin{definition}[Strong solution]\label{de2.1v3}
		A function $u: [0, T) \times \mathbb{R}^d  \to \mathbb{C}$ for $0< T \le \infty$ is a strong $\dot{H}^1$ solution of \eqref{eq1.1} if for all $t\in [0,T)$,
		$u \in \left(C_t^0 \dot{H}^1_x \cap L_{t,x}^\frac{2(d+2)}{d- 2} \right)([0,t] \times \mathbb{R}^d ) $, and obeys the Duhamel formula
		\begin{align*}
			u(t) = e^{t z \Delta } u_0 + \int_0^t e^{(t - s) z \Delta } z f(u(s)) \,\mathrm{d}s.
		\end{align*}
		
	\end{definition}
	For the linear part of the Ginzburg-Landau equation $\bar{z} \partial_t u = \Delta u $, it is well-known that the operator $z  \Delta $ with domain $H^2 \left( \mathbb{R}^d \right)$ generates a semi-group
	$ \left(e^{tz \Delta} \right)_{t\ge 0 }$ on $L^2 \left( \mathbb{R}^d \right)$; see for instance \cite{GJL,MYZ}. Moreover, the semi-group $ \left(e^{tz \Delta} \right)_{t \ge 0}$ is analytic for $\Re z > 0$ and
	\begin{align*}
		\left(e^{t z \Delta } \phi \right)(x) = \left(4 \pi  z t \right)^{- \frac{d}2} \left( e^{- \frac{ | \cdot |^2}{ 4 z t} } \ast \phi \right)(x).
	\end{align*}
	By the elementary estimates
	\begin{align*}
		\left\|  \left(4  \pi z t \right)^{- \frac{d}2}  e^{- \frac{ |x|^2}{ 4 z t} }  \right\|_{L^r_x}
		\lesssim
		\begin{cases}
			r^{- \frac{d}{2 \sigma }}
			( 4 \pi t)^{- \frac{d}2  \left( 1 - \frac1r \right) }  \left(\Re z \right)^{- \frac{d}{2 r } }, \ 1 \le r < \infty, \\
			( 4 \pi t)^{- \frac{d}2}, \ r = \infty,
		\end{cases}
	\end{align*}
	and Young's inequality, we infer that for $1 \le p \le r \le \infty$,
	\begin{align*}
		\left\| e^{t z \Delta }  \psi \right\|_{L_x^r } \lesssim  \left( \Re z \right)^{- \frac{d}2  \left( 1- \frac1p + \frac1r \right)} t^{- \frac{d}2
			\left( \frac1p - \frac1r \right) } \|\psi \|_{L_x^p}.
	\end{align*}
	Then, following the arguments in \cite{KT}, one easily obtains the  Strichartz estimates for Ginzburg-Landau equation, for which we formulate in below.
	
	\begin{lemma}
		\label{le2.2v21}
		The following Strichartz estimates hold.
		\begin{enumerate}
			{\upshape\item Homogeneous Strichartz estimate:}
			For $ 2\le q, r \le \infty$,  with $(q,r,d) \ne (2, \infty, 2)$ and $\frac2q = d\left( \frac12 - \frac1r \right)$, we have
			\begin{align}\label{eq2.2v3}
				\left\|e^{t z  \Delta} f \right\|_{L^q_t L^r_x}  \lesssim
				\|f\|_{L^2_x}.
			\end{align}
			{\upshape\item Inhomogeneous Strichartz estimate:}
			For $ 2\le q, r, \tilde{q},\tilde{r} \le \infty$,
			\begin{align}
				\left\|\int_0^t e^{(t - s) z \Delta } F(s) \,\mathrm{d}s \right\|_{L_t^q L_x^r ( \mathbb{R}_+ \times \mathbb{R}^d  )} \lesssim \| F \|_{L_t^{\tilde{q}'} L_x^{\tilde{r}'}}, \notag
			\end{align}
			where $\frac{2}q + \frac{d}r = \frac2{\tilde{q}} + \frac{d}{\tilde{r}} =  \frac{d}2$, with $(q,r,d) , \left(\tilde{q}, \tilde{r}, \tilde{d} \right) \ne (2, \infty, 2)$
		\end{enumerate}	
	\end{lemma}
	
	With the aid of Lemma \ref{le2.2v21}, the following well-posedness theory is immediate. Since the proof is rather standard, we do not recall it here and refer the interested readers to \cite{GJL,GV2,HW,W0} for details.
	
	\begin{theorem}[Local well-posedness]\label{th2.1v3}
		Assume $u_0 \in \dot{H}^1(\mathbb{R}^d )$. {{}Then the following conclusions hold:}
		\begin{enumerate}
			{\upshape\item Local existence:}
			There exists a unique, maximal-lifespan solution to \eqref{eq1.1} on $\left[0, T_{\max} (u_0) \right) \times \mathbb{R}^d $.
			
			{\upshape\item Blow-up criterion:}
			If $T_{\max} < \infty$, then $\|u \|_{S([0, T_{\max}(u_0)))} = \infty$.
			
			{\upshape\item Unconditional uniqueness:}
			If $u_1, u_2$ are two solutions of \eqref{eq1.1} on $[0, T)$ with $u_1 (0) = u_2(0) $, then $u_1 = u_2. $
			
			{\upshape\item Small data global existence:}
			There is $\epsilon_0 > 0$ such that if $\left\|e^{t z \Delta } u_0 \right\|_{S(\mathbb{R}_+)} \le \epsilon_0$, then the solution $u$ is global and
			$
			\|u \|_{S(\mathbb{R}_+)}
			\lesssim \epsilon_0.$
			This holds in particular when $\|u_0 \|_{\dot{H}^1}$ is sufficiently small.
		\end{enumerate}
		
	\end{theorem}
	
	The next stability result follows from the argument in \cite{CKSTT}.
	\begin{proposition}[Stability] \label{pr2.1v3}
		For any $E, L > 0$ and $\epsilon > 0$, there exists $\delta > 0$ such that if $\tilde{u}: [0, T) \times \mathbb{R}^d  \to \mathbb{C}$ is an approximate solution to \eqref{eq1.1} in the sense that
		\begin{align*}
			\|\nabla e\|_{L_{t,x}^\frac{2(d+2)}{d+4}  ( [0, T) \times \mathbb{R}^d )} \le \delta,
		\end{align*}
		where $e: = \bar{z}  \partial_t \tilde{u} - \Delta \tilde{u} -  f \left( \tilde{u}  \right) $	with
		\begin{align*}
			\left\|\tilde{u} \right\|_{L_t^\infty \dot{H}_x^1 ( [0, T) \times \mathbb{R}^d )} \le E
			\quad\text{ and }\quad
			\|\tilde{u} \|_{S([0, T) )} \le L,
		\end{align*}
		and if $u_0 \in \dot{H}^1( \mathbb{R}^d )$ satisfies
		\begin{align*}
			\|u_0 - \tilde{u}(0) \|_{\dot{H}^1} \le \delta ,
		\end{align*}
		then there exists a solution $u: [0, T) \times \mathbb{R}^d  \to \mathbb{C}$ of \eqref{eq1.1} with $u(0) = u_0$ such that
		\begin{align*}
			\| u - \tilde{u} \|_{L_t^\infty \dot{H}_x^1 } + \|u - \tilde{u} \|_{S([0, T)
				)} \le \epsilon.
		\end{align*}
		
	\end{proposition}
	
	\subsection{Some variational facts of Ginzburg-Landau equation}
	We have the following key facts on estimates of the derivatives of $E(u)$ and $\|u(t)\|_{L^2}^2$, which play an important role in the proofs of global well-posedness and finite time blow-up.
	\begin{proposition}\label{pr2.3v3'}
		Let $u_0 \in C^0 ( \mathbb{R}^d) \cap H^1( \mathbb{R}^d)$ and $u$ be the corresponding solution to \eqref{eq1.1} defined on the maximal interval $[0, T_{\max} )$.
		\begin{itemize}
			\item[ (i)]
			For any $0 \le s < t < T_{\max}$,
			\begin{align}\label{eq2.3v3'}
				E(u(s)) = E(u(t)) + \Re z \int_s^t \int_{ \mathbb{R}^d} |u_\tau|^2 \,\mathrm{d}x \mathrm{d} \tau,
			\end{align}
			which implies the energy $E(u(t))$ is a non-increasing function in time.
			\item[ (ii)]
			Set
			\begin{align*}
				K(u) = \int |\nabla u |^2 \,\mathrm{d}x  - \int |u |^{\frac{2d}{d-2} } \,\mathrm{d}x .
			\end{align*}
			For all $0 < t < T_{\max}$, we have
			\begin{align}\label{eq5.11v15}
				\frac{d}{dt} \int |u(t) |^2 \,\mathrm{d}x
				= - 2 \Re z K (u)
				= \frac4d \Re z \int  |u|^\frac{2d}{d-2} \,\mathrm{d}x - 4 \Re z E(u) .
			\end{align}
		\end{itemize}
		
	\end{proposition}
	\begin{proof}
		Multiplying \eqref{eq1.1} by $\bar{u}_t$, integrating by parts and taking the real part, we get
		\begin{align*}
			\Re z \int |u_t|^2 \,\mathrm{d}x
			=- \partial_t \left( \int \frac12 |\nabla u|^2 \,\mathrm{d}x - \mu \frac{d-2}{2d} \int |u|^\frac{2d}{d-2}  \,\mathrm{d}x \right).
		\end{align*}
		This yields
		\begin{align*}
			- \frac{d E(u(t))}{dt} =  \Re z \int |u_t|^2 \,\mathrm{d}x.
		\end{align*}
		By the fundamental theorem of calculus, we have \eqref{eq2.3v3'}.
		
		We now turn to \eqref{eq5.11v15}. Multiplying \eqref{eq1.1} by $\bar{z} \bar{u}$ and integrating by parts, we obtain
		\begin{align*}
			\Re \int_{\mathbb{R}^d} \bar{u} \partial_t u \,\mathrm{d}x = - \Re z  K(u) ,
		\end{align*}
		which is exactly \eqref{eq5.11v15}.
	\end{proof}
	
	The following liner profile decomposition describes the lack of compactness of the embedding $e^{ t z \Delta}: \dot{H}^1 \hookrightarrow  L_{t,x}^\frac{2(d+2)}{d-2}$,
	which follows from the argument of the linear profile decomposition of the Schr\"odinger equation \cite{Ke} and also the linear profile decomposition of the Navier-Stokes equation \cite{G,GKP}. This will play an important role in the proof of Theorem \ref{pr5.2v32}. To ease our notation, from now on, we write
	\begin{equation}\label{eq:scaled functions}
		\phi_{x_0,\lambda}(x):=\frac1{  \lambda^\frac{d-2 }2 } \phi \left( \frac{x - x_0}{ \lambda} \right)\quad \text{and}\quad v_{x_0,\lambda}(t,x):=\frac1{  \lambda^\frac{d-2 }2 } v \left(\frac{t}{\lambda^2}, \frac{x - x_0}{ \lambda} \right).
	\end{equation}
	
	\begin{proposition}[Linear profile decomposition in $\dot{H}^1$]\label{pr2.2v3}
		Let $\{u_n\}_n$ be a bounded sequence of functions in $\dot{H}^1 \left( \mathbb{R}^d \right)$. Then, after possibly passing to a subsequence, there exist $K\in \{ 1,2, \cdots, \infty\}$, functions $ \left\{\phi^j \right\}_{j= 1}^{K} \subseteq \dot{H}^1(\mathbb{R}^d)$ and $ \left\{ \left(  x_n^j , \lambda_n^j  \right) \right\} \subseteq  \mathbb{R}^d \times \mathbb{R}_+$ such that for any $0 \le k \le K$, we have the following decomposition
		\begin{align*}
			u_n (x) = \sum\limits_{j = 1}^k\phi^j_{x_n^j,\lambda_n^j}(x)+ w_n^k(x),
		\end{align*}
		with the following properties:
		\begin{itemize}
			\item[(i)] The error term satisfies
			\begin{align}
				\limsup\limits_{n \to \infty}
				\left\|e^{tz \Delta } w_n^k  \right\|_{S(\mathbb{R}_+ )	} \xrightarrow{k\to K} 0,  \label{eq4.4vnew}
				\intertext{ and }
				\left( \lambda_n^k  \right)^\frac{d-2}2 w_n^k   \left( \lambda_n^k  x + x_n^k  \right)
				\stackrel{n\to \infty}{\rightharpoonup} 0 \quad\text{ in } \dot{H}^1( \mathbb{R}^d). \notag
			\end{align}
			
			\item[(ii)] For all $j\neq l$, it holds
			\begin{align}\label{eq4.5vnew}
				\frac{ \lambda_n^j}{\lambda_n^l} + \frac{\lambda_n^l}{\lambda_n^j} + \frac{ \left|x_n^j -x_n^l \right| }{ \lambda_n^j  }
				\xrightarrow{n\to\infty} \infty.
			\end{align}

			\item[(iii)] For all $ 1 \le k < K  $, as $n \to \infty$, the following decoupling properties hold
			\begin{align}\label{eq2.9v50}
				\|u_n\|_{\dot{H}^1}^2 = \sum\limits_{j = 1}^k \left\|\phi^j \right\|_{\dot{H}^1}^2  + \left\|w_n^k \right\|_{\dot{H}^1}^2 + o_n(1),
			\end{align}
			and
			\begin{align*}
				E(u_n) = \sum\limits_{j = 1}^k E \left( \phi^j \right) + E \left(w_n^k \right) + o_n(1),
			\end{align*}
			where $o_n(1)\to 0$ as $n\to \infty$.
		\end{itemize}
	\end{proposition}
	
	Following the proofs of Theorem 3.9 and Remark 3.14 in \cite{KM}, one can prove the energy trapping.
	\begin{lemma}[Energy trapping]\label{le2.1v3}
		Fix a sufficiently small $\delta_0>0$. If
		\begin{align*}
			E(u_0) \le ( 1 - \delta_0) E(W), \
			\|\nabla u_0 \|_{L^2} <  \|\nabla W\|_{L^2},
		\end{align*}
		then there exists $\bar{\delta} = \bar{\delta}  ( \delta_0, d ) > 0$ such that for all $t \in [0, T_{\max} )$, the solution $u$ of \eqref{eq1.1} satisfies
		\begin{align}
			& \int |\nabla u(t)|^2 \,\mathrm{d}x \le  \left( 1- \bar{\delta } \right) \int |\nabla W|^2\, \mathrm{d}x, \notag  \\
			& K(u(t)) \ge \bar{\delta} \int |\nabla u(t)|^2 \,\mathrm{d}x, \label{eq2.11v3}
		\end{align}
		and $E(u(t)) \ge 0$.
		
		On the other hand, if $E(u_0) < E(W)$ and $ \|u_0 \|_{\dot{H}^1} > \|W\|_{\dot{H}^1}$, then
		\begin{align*}
			\|u(t) \|_{\dot{H}^1} > \|W\|_{\dot{H}^1},
		\end{align*}
		and there is a constant $\delta_3 > 0$ such that
		\begin{align}\label{eq2.19v41}
			K(u(t))
			\le - \delta_3.
		\end{align}				
	\end{lemma}

	\section{Proof of Blow-up }\label{se3v23}
	In this section, we shall adapt the arguments from \cite{Le} to prove the blow-up part of Theorem \ref{th1.1v3}.
	
	\begin{proof}[Proof of Theorem \ref{th1.1v3} (2)]
		Suppose for contradiction that $T_{\max}(u_0) = \infty$, and denote
		\begin{align*}
			I(t) : = \frac12 \int_0^t \|u(s) \|_{L_x^2}^2 \,\mathrm{d}s.
		\end{align*}
		Then $I'(t) = \frac12 \|u(t) \|_{L_x^2}^2$ {{}and it follows from \eqref{eq5.11v15} that}
		\begin{align*}
			I''(t) = - \Re z K(u) = \Re z \left( - \frac{2d }{d - 2} E(u) + \frac2{d-2 } \int |\nabla u |^2 \,\mathrm{d}x \right) .
		\end{align*}
		By \eqref{eq2.19v41}, we have
		\begin{align*}
			I''(t) \ge  \Re z \delta_3 > 0,
		\end{align*}
		and so $I'(t) \to \infty$ and $I(t) \to \infty$ as $t \to \infty$.
		
		Note that by \eqref{eq2.3v3'},
		\begin{align*}
			\Re z \int_0^t \|u_s (s) \|_{L_x^2}^2 \,\mathrm{d}s = E(u_0) - E(u(t)) \le  E(W)- E(u(t)),
		\end{align*}
		which implies
		\begin{align*}
			I''(t) &\ge - \Re z  \frac{2d}{d - 2} E(u)
			+ \Re z  \frac{2}{d- 2} \|\nabla W\|_{L^2}^2\\
			&= - \Re z  \frac{2d}{d - 2} E(u) + \Re z \frac{2d}{d- 2}E(W)
			\ge  ( \Re z )^2 \frac{2d }{d - 2} \int_0^t \|u_s (s) \|_{L_x^2}^2 \,\mathrm{d}s.
		\end{align*}
		Hence by H\"older's inequality and \eqref{eq5.11v15}, we have
		\begin{align*}
			I(t) I''(t) & \ge (\Re z )^2  \frac{d}{d- 2}  \left( \int_0^t \|u_s (s) \|_{L_x^2}^2 \,\mathrm{d}s  \right)  \left( \int_0^t \|u(s) \|_{L_x^2}^2 \,\mathrm{d}s  \right) \\
			& \ge ( \Re z )^2  \frac{d}{d - 2}  \left( \int_0^t \int  \left| \overline{ u(s, x)  } u_s (s, x) \right|  \,\mathrm{d}x \mathrm{d}s  \right)^2
			=   \frac{d}{d- 2}  \left( I'(t) - I'(0) \right)^2.
		\end{align*}
		Since $I'(t) \to \infty$ as $t \to \infty$, the last inequality implies existence of $\alpha$, $t_0> 0$, such that
		\begin{align*}
			I(t) I''(t) \ge ( 1 + \alpha)
			\left(I'(t) \right)^2, \text{ for } t \ge t_0.
		\end{align*}
		This inequality guarantees that the non-increasing function $I(t)^{- \alpha}$ is concave on $[t_0, \infty)$, which contradicts the fact that $I(t)^{- \alpha} \to 0, $ as $t \to \infty$.
	\end{proof}
	
	\section{Existence of a critical element}\label{se5v23}
	
	In this section, we show the existence of a critical element if the global existence part of Theorem \ref{th1.1v3} does not hold, following  the arguments in~\cite{KM,KM1,KK11,KV2010}. We show that if the critical energy of the initial data is strictly less than the energy of $W$ and if the kinetic energy of the initial data is strictly less than the kinetic energy of $W$, then there exists a critical element, which is pre-compact up to scaling and translation in $\dot{H}^1$.
	
	For any $0 \le E_0 \le E(W)$, we define
	\begin{align*}
		\Lambda(E_0) : = \sup \left\{ \|u\|_{S([0, T))}	:
		\text{ solutions $u$ of \eqref{eq1.1} on $[0,T) $ with $E(u(0) ) \le E_0$  and $\|u(0)\|_{\dot{H}^1} < \|W\|_{\dot{H}^1}$ } \right\},
	\end{align*}
	where $ [ 0, T)$ is  the maximal lifespan 
	of the solution $u$ and then set
	\begin{align*}
		E_c = \sup \left\{ E_0: \Lambda(E_0) < \infty \right\}.
	\end{align*}
	Suppose by contradiction that
	$E_c  < E(W)$, we can then take a sequence $\{u_n \}_{n \ge 1}$ of solutions (up to time translations) to \eqref{eq1.1} with maximal lifespan being $[0,T_n )$ such that
	\begin{align} \label{eq5.2v32}
		E(u_n(0) ) < E(W), E(u_n(0)) \xrightarrow{n \to \infty} E_c ,
		\|\nabla u_n(0) \|_{L^2} < \|\nabla W\|_{L^2} \text{ and } \|u_n \|_{S(0, T_n) } \xrightarrow{n \to \infty} \infty.
	\end{align}
	Applying Proposition \ref{pr2.2v3} to $\{ u_n(0) \}_{n \ge 1}$, {{}we obtain that}, passing to a subsequence if necessary, there exist $\phi^j \in \dot{H}^1$, $w_n^k \in \dot{H}^1$, and $ \left(\lambda_n^j,  x_n^j \right) \subseteq \mathbb{R}_+\times \mathbb{R}^d$ with the properties listed in Proposition \ref{pr2.2v3} such that
	\begin{align*}
		e^{tz \Delta} u_n(0)
		= \sum\limits_{j = 1}^k \phi_n^j(t,x) + e^{tz \Delta} w_n^k
		:= \sum\limits_{j = 1}^k e^{   t z \Delta}  \phi^j_{x_n^j,\lambda_n^j}(x)+ e^{tz \Delta} w_n^k.
	\end{align*}
	Let $v^j$ be the solution of \eqref{eq1.1} with $v^j(0,x) = \phi^j(x)$, where $I^j$ is the maximal lifespan of $v^j$. We also have
	\begin{align*}
		E \left( \phi^j
		\right) < E(W) \quad\text{ and }\quad \left\|
		\phi^j \right\|_{\dot{H}^1} < \|W\|_{\dot{H}^1}.
	\end{align*}
	Denote the nonlinear profile by
	\begin{align*}
		v_n^j(t,x) : = v^j_{x_n^j,\lambda_n^j}(t,x).
	\end{align*}
	Then we obtain the corresponding nonlinear profile decomposition
	\begin{align}\label{eq4.5newv44}
		u_n^k : = \sum\limits_{j = 1}^k v_n^j + e^{tz \Delta} w_n^k.
	\end{align}
	We next show that $u_n^k$ is an approximate solution of $u_n$ {{}for sufficiently large $n$, provided that each nonlinear profile has finite $S(\mathbb{R}_+)$-norm}.
	
	\begin{lemma}\label{le4.4v32}
		Let $v^j(t)$ be constructed as above,
		then there exists $j_0\in \mathbb{N}$ large enough such that
		$I^j = \mathbb{R}_+$ for $j \ge j_0$ and
		\begin{align*}
			\sum\limits_{j \ge j_0}  \left\|\nabla v^j \right\|_{L_t^\infty L_x^2 \cap L_t^2 L_x^\frac{2d}{d-2}   \left( \mathbb{R}_+ \times \mathbb{R}^d \right)}^2 \le \sum\limits_{j \ge j_0} \left\| \phi^j \right\|_{\dot{H}^1}^2 < \infty.
		\end{align*}
	\end{lemma}

	\begin{proof}
		For any $\delta > 0$, by \eqref{eq2.9v50}, there exist two constants $j_0$ and $n_0$ large enough such that if $n \ge n_0$,
		\begin{align*}
			\sum\limits_{j \ge j_0}  \left\|  \phi^j  \right\|_{\dot{H}^1}^2 < \delta.
		\end{align*}
		By Lemma \ref{le2.2v21}, Leibnitz's rule, H\"older's inequality and Sobolev's inequality, we have
		\begin{align*}
			\left\| \nabla v^j \right\|_{L_t^\infty L_x^2 \cap L_t^2 L_x^\frac{2d}{d-2}  \left( \mathbb{R}_+ \times \mathbb{R}^d  \right)}
			\lesssim  \left\|\phi^j \right\|_{\dot{H}^1 } +  \left\|\nabla v^j  \right\|_{L_t^\infty L_x^2 \cap L_t^2 L_x^\frac{2d}{d-2 } }^\frac{d+2}{d-2} .
		\end{align*}
		Hence, by the continuity argument, for sufficiently small $\delta > 0$, we have
		\begin{align*}
			\left\| \nabla v^j \right\|_{L_t^\infty L_x^2 \cap L_t^2 L_x^\frac{2d}{d-2} } \lesssim  \left\| \phi^j \right\|_{\dot{H}^1}.
		\end{align*}
		Thus,
		\begin{align*}
			\sum\limits_{j \ge j_0} \left\|\nabla v^j \right\|_{L_t^\infty L_x^2 \cap L_t^2 L_x^\frac{2d}{d- 2} }^2 \lesssim \sum\limits_{j \ge j_0} \left\|\phi^j \right\|_{\dot{H}^1}^2.
		\end{align*}
		
	\end{proof}
	We also have
	\begin{lemma}\label{le4.5v32}
		Fix $k_0 \in \mathbb{N}$ and suppose for all $1 \le j \le k_0$,
		\begin{align}\label{eq4.24v32}
			\left\|v^j \right\|_{L_{t,x}^\frac{2(d+2)}{ d- 2}  \left(I^j \times \mathbb{R}^d  \right) } < \infty.
		\end{align}
		Then $I^j = \mathbb{R}_+$ and for $1 \le j \le k_0$, we have
		\begin{align}\label{eq4.25v32}
			\left\| \nabla v_n^j \right\|_{L_{t,x}^\frac{2(d+2)}d ( \mathbb{R}_+ \times \mathbb{R}^d )} + \left\|v_n^j \right\|_{L_{t,x}^\frac{2(d + 2)}{ d - 2}  ( \mathbb{R}_+ \times \mathbb{R}^d )} \lesssim \left\|\nabla v^j \right\|_{L_t^\infty L_x^2 \cap L_t^2 L_x^\frac{2d}{d- 2}  ( \mathbb{R}_+ \times \mathbb{R}^d )} \lesssim 1.
		\end{align}
		Moreover, for each $k \in \mathbb{N}$ {{}with} $1 \le k \le k_0$, there exists $N_k \in \mathbb{N}$ such that
		\begin{align}\label{eq4.26v32}
			\sup\limits_{n \ge N_k} \Bigg(  \left\|  \nabla \Big(
			\sum\limits_{j = 1}^k   v_n^j  \Big)
			\right\|_{L_{t,x}^\frac{2(d+2)}d ( \mathbb{R}_+ \times \mathbb{R}^d )} +  \left\|\sum\limits_{j = 1}^k v_n^j \right\|_{L_{t,x}^\frac{2(d+2)}{d - 2}   \left( \mathbb{R}_+ \times \mathbb{R}^d  \right) } \Bigg)
			\lesssim 1 .
		\end{align}
		Furthermore, if \eqref{eq4.24v32} holds for {{}all} $j \in \mathbb{N}$, then
		\begin{align}\label{eq4.27v32}
			\limsup\limits_{n \to \infty}  \left\|v_n^j \nabla e^{tz \Delta} w_n^k  \right\|_{L_{t,x}^\frac{d+2}{d- 1} } \xrightarrow{k \to \infty} 0.
		\end{align}
		
	\end{lemma}
	\begin{proof}
		By Theorem \ref{th2.1v3} and \eqref{eq4.24v32}, we have	$I^j = \mathbb{R}_+$. Next, for the first inequality in \eqref{eq4.25v32}, {{}note that}
		\begin{align*}
			\left\| v_n^j  \right\|_{L_{t,x}^\frac{2(d+2)}{d - 2}  ( \mathbb{R}_+ \times \mathbb{R}^d )} =  \left\|v^j  \right\|_{L_{t,x}^\frac{2(d+2)}{d - 2} ( \mathbb{R}_+ \times \mathbb{R}^d )} \lesssim  \left\|\nabla v^j \right\|_{L_t^\infty L_x^2 \cap L_t^2 L_x^\frac{2d}{d - 2}  ( \mathbb{R}_+ \times \mathbb{R}^d )},
		\end{align*}
		and
		\begin{align*}
			\left\| \nabla v_n^j  \right\|_{L_{t,x}^\frac{2(d+2)}d ( \mathbb{R}_+ \times \mathbb{R}^d )} =  \left\|\nabla v^j  \right\|_{L_{t,x}^\frac{2(d+2)}d ( \mathbb{R}_+ \times \mathbb{R}^d )} \lesssim  \left\|\nabla v^j  \right\|_{L_t^\infty L_x^2 \cap L_t^2 L_x^\frac{2d}{ d - 2}(\mathbb{R}_+ \times \mathbb{R}^d ) }.
		\end{align*}
		Hence, to deduce \eqref{eq4.25v32}, it suffices to show
		\begin{align*}
			\left\|\nabla v^j  \right\|_{L_t^\infty L_x^2 \cap L_t^2 L_x^\frac{2d}{ d - 2} ( \mathbb{R}_+ \times \mathbb{R}^d )} \lesssim 1.
		\end{align*}
		By the Strichartz estimate, H\"older's inequality and \eqref{eq4.24v32}, we have
		\begin{align*}
			\left\| \nabla v^j \right\|_{L_t^\infty L_x^2 \cap L_t^2 L_x^\frac{2d}{d- 2}  (\mathbb{R}_+\times \mathbb{R}^d)}
			\lesssim  \left\|\phi^j  \right\|_{\dot{H}^1} +  
			\left\|\nabla v^j  \right\|_{L_t^\infty L_x^2 \cap L_t^2 L_x^\frac{2d}{d-2}}^\frac{d+2}{d - 2} .
		\end{align*}
		In this case, similar to the proof of Lemma \ref{le4.4v32}, we infer from the fact $ \left\|\phi^j  \right\|_{\dot{H}^1} \lesssim 1$ that
		\begin{align*}
			\left\|\nabla v^j  \right\|_{L_t^\infty L_x^2( \mathbb{R}_+ \times \mathbb{R}^d )} \lesssim 1,
		\end{align*}
		from which we conclude $ \left\|\nabla v^j  \right\|_{L_t^2 L_x^\frac{2d}{d - 2}  \cap L_t^\infty L_x^2 ( \mathbb{R}_+ \times \mathbb{R}^d )} \lesssim 1$.
		
		To prove \eqref{eq4.26v32}, we apply the following elementary inequality
		\begin{align}\label{eq4.7v70}
			\left|  \left|\sum\limits_{j = 1}^k v_n^j \right|^\frac{2(d+2)}{d - 2}  - \sum\limits_{j = 1}^k  \left|v_n^j \right|^\frac{2(d+2)}{d - 2}   \right|
			\lesssim_k \sum\limits_{ j = 1}^k \sum\limits_{\substack{1 \le j' \le k \\ j' \ne j}}  \left|v_n^j \right|^\frac{d+6}{d - 2}   \left|v_n^{j'} \right|,
		\end{align}
		to obtain
		\begin{align}\label{eq4.34v32}
			\left\| \sum\limits_{j = 1}^k v_n^j  \right\|_{L_{t,x}^\frac{2(d+2)}{d - 2} ( \mathbb{R}_+ \times \mathbb{R}^d )}^\frac{2(d+2)}{d - 2}
			\le \sum\limits_{j = 1}^k  \left\|v_n^j  \right\|_{L_{t,x}^\frac{2(d+2)}{d - 2} }^\frac{2(d+2)}{d - 2}
			+ C_k \sum\limits_{j = 1}^k \sum\limits_{ \substack{ 1 \le j' \le k \\  j \ne j'}} \int_{\mathbb{R}_+} \int_{\mathbb{R}^d }  \left|v_n^j \right|^\frac{d+6}{d - 2}   \left|v_n^{j'} \right| \,\mathrm{d}x \mathrm{d}t.
		\end{align}
		Without loss of generality, we assume that $k \ge j_0$ with $j_0$ given by Lemma \ref{le4.4v32}. Then, we have
		\begin{equation}\label{eq4.35v32}
			\begin{aligned}
				\sum\limits_{j = 1}^k  \left\|v_n^j \right\|_{L_{t,x}^\frac{2(d+2)}{d - 2}  ( \mathbb{R}_+ \times \mathbb{R}^d )}^\frac{2(d+2)}{d - 2}
				\lesssim \sum\limits_{1 \le j \le j_0}  \left\|\nabla v^j \right\|_{L_t^\infty L_x^2 \cap L_t^2 L_x^\frac{2d}{d - 2} }^\frac{2(d+2)}{d -2 }
				+ \sum\limits_{j_0 < j \le k} \left\|\phi^j \right\|_{\dot{H}^1}^2 < \infty.
			\end{aligned}
		\end{equation}
		On the other hand, by \eqref{eq4.5vnew}, there exists $N_k \in \mathbb{N}$ such that for $n \ge N_k$ and all $j \ne j'$
		\begin{align}\label{eq4.36v32}
			\int_{\mathbb{R}_+} \int_{\mathbb{R}^d }  \left|v_n^j \right|^\frac{d+6}{d - 2}   \left|v_n^{j'} \right| \,\mathrm{d}x \mathrm{d}t \le \frac1{C_k k^2}.
		\end{align}
		Hence, by \eqref{eq4.34v32}, \eqref{eq4.35v32}, and \eqref{eq4.36v32}, we deduce that for any $1 \le k \le k_0$,
		there exists an $N_k' \in \mathbb{N}$ such that
		\begin{align*}
			\sup\limits_{n \ge N_k'}  \left\|\sum\limits_{j = 1}^k v_n^j  \right\|_{L_{t,x}^\frac{2(d+2)}{d - 2}  ( \mathbb{R}_+ \times \mathbb{R}^d )} \lesssim 1 .
		\end{align*}
		Similarly, applying  \eqref{eq4.7v70}, Lemma \ref{le4.4v32} and \eqref{eq4.5vnew},
		we infer that for any $1 \le k \le k_0$, there is an $ N_k^{''} \in \mathbb{N}$ such that
		\begin{align*}
			\sup\limits_{n \ge N_k^{''}}  \left\| \sum\limits_{j = 1}^k \nabla v_n^j  \right\|_{L_{t,x}^\frac{2(d+2)}d ( \mathbb{R}_+ \times \mathbb{R}^d )} \lesssim 1.
		\end{align*}
		Hence, \eqref{eq4.26v32} follows by taking $N_k = \max \left( N_k^{'}, N_k^{''} \right)$.
		
		Finally, we turn to the proof of \eqref{eq4.27v32}.
		For any $\epsilon > 0$, we can find $\tilde{v}^j \in C_c^\infty( \mathbb{R}_+ \times \mathbb{R}^d )$ and $N_j \in \mathbb{N}$ such that
		\begin{align*}
			\left\|v^j - \tilde{v}^j  \right\|_{L_{t,x}^\frac{2(d+2)}{d -2 } ( \mathbb{R}_+ \times \mathbb{R}^d)} \le \epsilon
		\end{align*}
		for all $j \ge N_j$. Hence, by H\"older's inequality and the Strichartz estimate, we have
		\begin{equation}\label{eq4.40v32}
			\begin{aligned}
				\left\| v_n^j \nabla e^{tz \Delta} w_n^k  \right\|_{L_{t,x}^\frac{d+2}{d - 1} }
				\lesssim &   \left\|v^j  \left( \nabla e^{\cdot z \Delta} w_n^k \right)  \left(  \left( \lambda_n^j \right)^2 t, \lambda_n^j x + x_n^j \right)    \left( \lambda_n^j \right)^\frac{d}2   \right\|_{L_{t,x}^\frac{d+2}{d - 1} } \\
				\lesssim  &  \left\|v^j - \tilde{v}^j  \right\|_{L_{t,x}^\frac{2(d+2)}{d - 2} } \left\|w_n^k  \right\|_{\dot{H}^1} +  \left( \lambda_n^j \right)^\frac{d}2 \left\|\tilde{v}^j \left( \nabla e^{\cdot z \Delta} w_n^k  \right) \left( \left( \lambda_n^j \right)^2 t , \lambda_n^j x + x_n^j  \right) \right\|_{L_{t,x}^\frac{d+2}{d - 1} }.
			\end{aligned}
		\end{equation}
		Next, we assume $\supp(\tilde{v}^j) \subseteq  \left\{(t,x) \in \mathbb{R}_+ \times \mathbb{R}^d : 0 \le t \le T_j ,  |x| \le R_j  \right\}$ and set
		\begin{align*}
			\tilde{w}_n^{j,k}(x) : = w_n^k  \left( \lambda_n^j x + x_n^j \right).
		\end{align*}
		To estimate the second term on the right hand side of \eqref{eq4.40v32}, by H\"older's inequality, it suffices to show
		\begin{align*}
			\limsup\limits_{n \to \infty} \left( \lambda_n^j \right)^\frac{d}2 \left\| \left( \nabla e^{ \left( \lambda_n^j \right)^2 t	z \Delta} w_n^k \right) \left( \lambda_n^j x + x_n^j \right) \right\|_{L_{t,x}^2  \left( 0 \le t \le T_j, |x | \le R_j \right)} \xrightarrow{k \to \infty} 0.
		\end{align*}
		Note that by Proposition \ref{pr2.2v3}, we have
		\begin{equation}\label{eq4.41v32}
			\begin{aligned}
				&   \limsup\limits_{n \to \infty } \left( \lambda_n^j \right)^\frac{d}2
				\left\| \left( \nabla e^{ \left( \lambda_n^j \right)^2 t  z \Delta} w_n^k \right)
				\left( \lambda_n^j x + x_n^j \right) \right\|_{L_{t,x}^2   \left( 0 \le t \le T_j, |x | \le R_j \right)}\\
				= &  \limsup\limits_{n \to \infty }
				\left( \lambda_n^j  \right)^{\frac{d}2 - 1} \left\|\nabla e^{tz \Delta} \tilde{w}_n^{j,k}  \right\|_{L_{t,x}^2  ( 0 \le t \le T_j, |x| \le R_j )} \\
				\lesssim &	\limsup\limits_{n \to \infty }
				( \lambda_n^j)^{\frac{d}2 - 1}  T_j^\frac{2}{ 3(d+2)}  R_j^\frac{3d+2}{ 6(d+2)}  \left\|e^{tz\Delta} \tilde{w}_n^{j,k}  \right\|_{L_{t,x}^\frac{2(d+2)}{d - 2} }^\frac13  \left\|\nabla \tilde{w}_n^{j,k} \right\|_{L_{ x}^2}^\frac23\\
				\lesssim &
				\limsup\limits_{n \to \infty }
				T_j^\frac{2}{ 3(d+2)}  R_j^\frac{3d+2}{ 6(d+2)}  \left\|e^{tz \Delta} w_n^k  \right\|_{L_{t,x}^\frac{2(d+2)}{d - 2} }^\frac13  \left\|\nabla w_n^k  \right\|_{L_x^2}^\frac23 \xrightarrow{k \to \infty} 0.
			\end{aligned}
		\end{equation}
		Thus, \eqref{eq4.27v32} follows from \eqref{eq4.40v32} and \eqref{eq4.41v32}.
	\end{proof}
	
	By stability, we can prove	
	\begin{proposition}\label{pr4.6v32}
		{{}If for all $j\geq 1$},
		\begin{align}\label{eq4.43v32}
			\left\|v^j  \right\|_{ S(I^j)} < \infty,
		\end{align}
		then {{}for all $n\in \mathbb{N}$ sufficiently large, we have}
		\begin{align}\label{eq4.44v32}
			\left\|u_n \right\|_{ S(\mathbb{R}_+)} < \infty.
		\end{align}
		
	\end{proposition}
	
	\begin{proof}
		We shall prove that $u_n^k$ is a sequence of approximation solutions of $u_n$.
		First of all, by Lemma \ref{le4.5v32} and \eqref{eq4.43v32},
		we know that $v_n^j$ is a global solution and thus $u_n^k$ is global.
		Next, by \eqref{eq4.26v32} and $\dot{H}^1$-boundedness of $w_n^k$, we may find a constant $B >0$ with the following property:
		there exists $N_{1,k} > 0$ such that
		\begin{align*}
			\sup\limits_{n \ge N_{1,k}} \left\|u_n^k  \right\|_{L_{t,x}^\frac{2(d+2)}{d - 2} ( \mathbb{R}_+ \times \mathbb{R}^d  )} \le B.
		\end{align*}
		In fact, by Proposition \ref{pr2.2v3}, we have
		\begin{align*}
			\left\|u_n(0) - u_n^k (0)  \right\|_{\dot{H}^1} \le \sum\limits_{j = 1}^k  \left\|v_n^j(0) - \phi_n^j(0)  \right\|_{\dot{H}^1}.
		\end{align*}
		This implies that for any $k \in \mathbb{N}$, there exists $N_{2,k} > 0$ such that
		\begin{align*}
			\sup\limits_{n \ge N_{2,k}}  \left\|u_n(0) - u_n^k(0) \right\|_{\dot{H}^1} \le 1.
		\end{align*}
		Hence, if we were able to derive that there exist $k_0, N_0 \in \mathbb{N}$ such that when $k\ge k_0$ and $n\ge N_0$, \begin{align*}
			\left\|\nabla  \left( \left( \bar{z} \partial_t - \Delta \right) u_n^k -   f\left( u_n^k \right) \right) \right\|_{L_{t,x}^\frac{2(d+2)}{d + 4}  } \le \delta,
		\end{align*}
		then \eqref{eq4.44v32} would follow from Proposition \ref{pr2.1v3}.
		
		To verify the preceding assertion, we first note that
		\begin{align*}
			\left( \bar{z} \partial_t - \Delta \right) u_n^k -  f \left( u_n^k \right)
			= - f \left( u_n^k \right) + f \left( \sum\limits_{j = 1}^k v_n^j \right)
			+ \sum\limits_{j = 1}^k  \left( \bar{z} \partial_t - \Delta  \right) v_n^j -  f\left( \sum\limits_{j = 1}^k v_n^j \right).
		\end{align*}
		Hence, it suffices to show
		\begin{align}\label{eq4.51v32}
			\lim\limits_{k \to \infty} \limsup\limits_{n \to \infty} \left\|\nabla  \left( f \left( u_n^k \right) - f \left( \sum\limits_{j = 1}^k v_n^j  \right) \right) \right\|_{L_{t,x}^\frac{2(d+2)}{d+4 } } = 0,
		\end{align}
		and
		\begin{align}\label{eq4.52v32}
			\lim\limits_{k \to \infty} \limsup\limits_{n\to \infty}  \left\|\nabla  \left( \sum\limits_{j = 1}^k \left( \bar{z} \partial_t - \Delta  \right) v_n^j - f \left( \sum\limits_{j = 1}^k v_n^j  \right) \right)  \right\|_{L_{t,x}^\frac{2(d+2)}{d + 4} } = 0.
		\end{align}
		To establish \eqref{eq4.51v32}, we estimate, using H\"older's inequality, as follows.
		\begin{align*}
			&  \left\|\nabla \left( f \left(u_n^k \right) - f \left( \sum\limits_{j = 1}^k v_n^j  \right) \right)  \right\|_{L_{t,x}^\frac{2(d+2)}{d + 4} }
			=   \left\|\nabla  \left( f \left( \sum\limits_{j = 1}^k v_n^j + e^{tz \Delta } w_n^k  \right) - f \left( \sum\limits_{j = 1}^k v_n^j  \right) \right) \right\|_{L_{t,x}^\frac{2(d+2)}{d + 4} }\\
			\lesssim & \left\| \nabla  e^{tz \Delta } w_n^k  \right\|_{L_{t,x}^\frac{2(d+2)}d }  \left\|  e^{tz \Delta } w_n^k \right\|_{L_{t,x}^\frac{2(d+2)}{d - 2}}^\frac4{d-2}
			+  \left\| \nabla \sum\limits_{j = 1}^k v_n^j  \right\|_{L_{t,x}^\frac{2(d+2)}d }   \left\|  e^{tz \Delta } w_n^k  \right\|_{L_{t,x}^\frac{2(d+2)}{d- 2}}^\frac4{d-2}\\
			&+  \left\|  \left| \sum\limits_{j = 1}^k v_n^j \right|^\frac4{d-2} \nabla e^{tz \Delta } w_n^k  \right\|_{L_{t,x}^\frac{2(d+2)}{d + 4}}
			+  \left\|\nabla \sum\limits_{j = 1}^k v_n^j  \right\|_{L_{t,x}^\frac{2(d+2)}d}  \left\| \sum\limits_{j = 1}^k v_n^j  \right\|_{L_{t,x}^\frac{2(d+2)}{d - 2}}^\frac{6 - d }{d - 2}  \left\| e^{tz \Delta } w_n^k  \right\|_{L_{t,x}^\frac{2(d+2)}{d - 2}} \\
			=:  &  A_{1,n, k}  + A_{2, n, k}  + A_{3, n, k}  + A_{4, n, k} .
		\end{align*}
		Notice that when $d \le 5$, by \eqref{eq4.26v32} and \eqref{eq4.4vnew}, we have $\limsup\limits_{n \to \infty} A_{4,n, k} \xrightarrow{k\to \infty} 0$, and,
		when $d \ge 6$, the term $A_{4,n, k } $ disappears.
		
		We first deal with the terms $A_{1,n, k}$ and $A_{2,n,k}$. By H\"older's inequality, we have
		\begin{align}\label{eq4.17vnew}
			\left\|  \left|e^{tz \Delta} w_n^k \right|^\frac4{d-2} \left|\nabla e^{tz \Delta} w_n^k \right|  \right\|_{L_{t,x}^\frac{2(d+2)}{d + 4}}
			\lesssim  \left\|e^{tz \Delta} w_n^k  \right\|_{L_{t,x}^\frac{2(d+2)}{d-2} }^\frac4{d-2}
			\left\| \nabla e^{tz \Delta} w_n^k  \right\|_{L_{t,x}^\frac{2(d+2)}d }.
		\end{align}
		Then, by the Strichartz estimate, \eqref{eq4.17vnew},
		and \eqref{eq4.4vnew}, we have
		\begin{align*}
			\limsup\limits_{n \to \infty} A_{1,n, k }  \lesssim  \limsup\limits_{n \to \infty}  \left\|e^{tz \Delta} w_n^k  \right\|_{L_{t,x}^\frac{2(d+2)}{d - 2} }^\frac4{d-2}  \xrightarrow{k\to\infty} 0.
		\end{align*}
		Similarly, by \eqref{eq4.26v32}, one can prove that
		\begin{align*}
			\sup\limits_{n \ge N_k}  \left\| \nabla \left( \sum\limits_{j = 1}^k v_n^j \right)  \right\|_{L_{t,x}^\frac{2(d+2)}d ( \mathbb{R}_+ \times \mathbb{R}^d )} \lesssim 1,
		\end{align*}
		which provides a desired control on $A_{2,n, k } $.
		
		We next turn to the term $A_{3, n, k } $. By H\"older's inequality and \eqref{eq4.26v32}, we have
		\begin{align*}
			\left\| \left|\sum\limits_{j = 1}^k v_n^j \right|^\frac4{d-2}   \left|\nabla e^{tz \Delta} w_n^k \right|  \right\|_{L_{t,x}^\frac{2(d+2)}{d + 4} }
			&\lesssim  \left\|\sum\limits_{j = 1}^k v_n^j  \right\|_{L_{t,x}^\frac{2(d+2)}{d -2} }^\frac3{d-2}\cdot  \Big\|\nabla e^{tz \Delta} w_n^k \Big\|_{L_{t,x}^\frac{2(d+2)}d}^\frac{d-3}{d-2}\cdot   \left\|\sum\limits_{j = 1}^k v_n^j \nabla e^{tz \Delta} w_n^k  \right\|_{L_{t,x}^\frac{d+2}{d - 1}}^\frac1{d-2}\\
			&\lesssim  \left\|\sum\limits_{j = 1}^k v_n^j  \nabla e^{tz \Delta } w_n^k  \right\|_{L_{t,x}^\frac{d+2}{d - 1}}^\frac1{d-2}.
		\end{align*}
		Hence, to establish \eqref{eq4.51v32}, we just need to show
		\begin{align}\label{eq4.66v32}
			\lim\limits_{k \to \infty} \limsup\limits_{n \to \infty}  \left\|\sum\limits_{j = 1}^k v_n^j \nabla e^{tz \Delta} w_n^k  \right\|_{L_{t,x}^\frac{d+2}{d - 1} } = 0.
		\end{align}
		By Lemma \ref{le4.4v32}, for any $\epsilon > 0$, there exists $J(\epsilon) > 0$ such that
		\begin{align*}
			\sum\limits_{J(\epsilon) \le j \le k }  \left\| v_n^j  \right\|_{L_{t,x}^\frac{2(d+2)}{d - 2}  ( \mathbb{R}_+  \times \mathbb{R}^d )}^2
			\le \sum\limits_{j \ge J(\epsilon)}  \left\|\nabla v^j  \right\|_{L_t^\infty L_x^2 \cap L_t^2 L_x^\frac{2d}{d - 2} }^2 < \frac{\epsilon^2}4.
		\end{align*}
		By a similar argument as used in the proof of \eqref{eq4.26v32} and by \eqref{eq4.5vnew},
		there is a constant $N_{k, \epsilon } > 0$ such that when $n > N_{k, \epsilon}$, we have
		\begin{align*}
			\left\| \sum\limits_{J(\epsilon) \le j \le k } v_n^j  \right\|_{L_{t,x}^\frac{2(d+2)}{d -2 }  ( \mathbb{R}_+ \times \mathbb{R}^d )}
			\le \left( \sum\limits_{J(\epsilon ) \le j \le k } \|v_n^j \|_{L_{t,x}^\frac{2(d+2)}{d - 2}  ( \mathbb{R}_+ \times \mathbb{R}^d )}^2  \right)^\frac12 + \frac\epsilon2 < \epsilon,
		\end{align*}
		where in the last step, we used the boundedness of $ \left\|\phi^j  \right\|_{\dot{H}^1}$. Then, by H\"older's inequality and the Strichartz estimate, we infer that
		\begin{align}\label{eq4.69v32}
			\left\| \sum\limits_{J( \epsilon ) \le j \le k } v_n^j \nabla e^{tz \Delta} w_n^k  \right\|_{L_{t,x}^\frac{d+2}{d - 1}  ( \mathbb{R}_+  \times \mathbb{R}^d )}
			\le  \left\|\sum\limits_{J( \epsilon) \le j \le k } v_n^j  \right\|_{L_{t,x}^\frac{2(d+2)}{d - 2}  }  \left\|\nabla e^{tz \Delta} w_n^k  \right\|_{L_{t,x}^\frac{2(d+2)}d  } \lesssim \epsilon .
		\end{align}
		On the other hand, by \eqref{eq4.27v32}, we have
		\begin{align*}
			\lim\limits_{k \to \infty} \limsup\limits_{n \to \infty}  \left\|\sum\limits_{1 \le j \le J(\epsilon)} v_n^j \nabla e^{tz \Delta} w_n^k  \right\|_{L_{t,x}^\frac{d+2}{d - 1} } = 0,
		\end{align*}
		which together with \eqref{eq4.69v32} implies \eqref{eq4.66v32}. This concludes the proof of \eqref{eq4.51v32}.
		
		We now turn to the proof of \eqref{eq4.52v32}. Note that
	\begin{align}
	& \left\|\nabla  \left( \sum\limits_{j = 1}^k \left( \bar{z} \partial_t  -  \Delta  \right) v_n^j  -  f \left( 
		u_n^{k, app } 
		\right) \right)  \right\|_{L_{t,x}^\frac{2(d+2)}{d + 4} }  \nonumber
	\\
	& \lesssim    \left\|\nabla  \left( \sum\limits_{j = 1}^k f \left(v_n^j \right) - f \left( 
	u_n^{k, app } 
	\right) \right) \right\|_{L_{t,x}^\frac{2(d+2)}{d + 4}}
	\lesssim    \sum\limits_{\substack{ 1 \le j, j' \le k , \\ j \ne j'}}   \left\|  \left|v_n^j  \right|^\frac4{d-2}  \nabla v_n^{j'}  \right\|_{L_{t,x}^\frac{2(d+2)}{d + 4} }. \label{eq4.74v32}
\end{align}
		To estimate the right hand side of \eqref{eq4.74v32}, we use density. 
		For any $\epsilon > 0$, there exists $\tilde{U}^j, \tilde{V}^j \in C_c^\infty( \mathbb{R}_+  \times \mathbb{R}^d )$ such that for sufficiently large $n$,
\begin{align*}
	\left\| v^j_{x_n^j,\lambda_n^j}(t,x)- \tilde{U}^j_{x_n^j,\lambda_n^j}(t,x)\right\|_{ S( \mathbb{R}_+ ) 
	}  +  \left\| \nabla   \left( v^j_{x_n^j,\lambda_n^j}(t,x)-\tilde{V}_{x_n^j,\lambda_n^j}^j(t,x) \right)   \right\|_{L_{t,x}^\frac{2(d+2)}d   \left( \mathbb{R}_+  \times \mathbb{R}^d  \right)} < \epsilon.
\end{align*}
		Hence, the right hand side of \eqref{eq4.74v32} is bounded from above by
		\begin{align*}
			\epsilon + \sum\limits_{ \substack{ 1 \le j, j' \le k\\ j \ne j'}}
			\left\|  \left| \tilde{U}^j_{x_n^j,\lambda_n^j}(t,x) \right|^\frac4{d-2}  \nabla \tilde{V}^{j'}_{x_n^{j'},\lambda_n^{j'}}(t,x) \right\|_{L_{t,x}^\frac{2(d+2)}{d + 4} }.
		\end{align*}
		If $\frac{ \lambda_n^j}{ \lambda_n^{j'}} \to \infty$, {{}then} we have
		\begin{align*}
			& \left\|  \left| \tilde{U}^j_{x_n^j,\lambda_n^j}(t,x) \right|^\frac4{d-2}  \nabla \tilde{V}^{j'}_{x_n^{j'},\lambda_n^{j'}}(t,x) \right\|_{L_{t,x}^\frac{2(d+2)}{d + 4} }
			\\
			\lesssim & \left( \lambda_n^{j'}  \right)^2  \left( \lambda_n^j \right)^{- 2}
			\left\|   \left| \tilde{U}^j \left( \frac{  \left( \lambda_n^{j'}  \right)^2 t }{  \left( \lambda_n^j \right)^2 }, \frac{ \lambda_n^{j'} x + x_n^{j'} - x_n^j }{ \lambda_n^j } \right) \right|^\frac4{d-2}  \left( \nabla \tilde{V}^{j'}  \right)(t,x) \right\|_{L_{t,x}^\frac{2(d+2)}{d+4} } \notag \\
			\lesssim  & \left( \lambda_n^{j'} \right)^2  \left( \lambda_n^j \right)^{- 2}  \left\| \tilde{U}^j \right\|_{L_{t,x}^\infty}^\frac4{d-2}  \left\|\nabla \tilde{V}^{j'} \right\|_{L_{t,x}^\frac{2(d+2)}{d + 4} }  \xrightarrow{n\to\infty} 0. \notag
		\end{align*}
		If $\frac{ \lambda_n^j}{\lambda_n^{j'}} \to 0$, {{}then} we have
		\begin{align*}
			& \left\| \left| \tilde{U}^j_{x_n^j,\lambda_n^j}(t,x) \right|^\frac4{d-2 }  \nabla \tilde{V}^{j'}_{x_n^{j'},\lambda_n^{j'}}(t,x)  \right\|_{L_{t,x}^\frac{2(d+2)}{d + 4} } \\
			\lesssim &   \left( \lambda_n^j \right)^2  \left( \lambda_n^{j'}  \right)^{-2}  \left\|  \left|\tilde{U}^j ( t,x) \right|^\frac4{d-2}   \left( \nabla \tilde{V}^{j'}  \right) \left( \frac{  \left( \lambda_n^j \right)^2 t }{  \left( \lambda_n^{j'} \right)^2} , \frac{  \lambda_n^jx + x_n^j -x_n^{j'} }{ \lambda_n^{j'} }  \right)  \right\|_{L_{t,x}^\frac{2(d+2)}{ d + 4} }
			\lesssim     \left( \lambda_n^j \right)^2  \left( \lambda_n^{j'} \right)^{-2 }  \xrightarrow{n\to\infty} 0.
		\end{align*}
		If there exists some $C_j > 0$ such that for any $n$,
		\begin{align*}
			\frac{ \lambda_n^j}{ \lambda_n^{j'}} + \frac{ \lambda_n^{j'}}{ \lambda_n^j} \le C_j,
		\end{align*}
		then we have by \eqref{eq4.5vnew} that
		\begin{align*}
			\left|\frac{ x_n^j - x_n^{j'} }{ \lambda_n^j}  \right|  \xrightarrow{n\to\infty} \infty.
		\end{align*}
		By symmetry, we may assume that $j < j'$. For sufficiently large $n$, we have
		\begin{align*}
			&  \left\|  \left|\tilde{U}^j_{x_n^j,\lambda_n^j}(t,x) \right|^\frac4{d-2}  \nabla \tilde{V}^{j'}_{x_n^{j'},\lambda_n^{j'}}(t,x)   \right\|_{L_{t,x}^\frac{2(d+2)}{d + 4}  } \\
			\lesssim  &  \left( \lambda_n^{j'} \right)^2  \left( \lambda_n^j \right)^{-2}  \left\|  \left|\tilde{U}^j \left( \frac{  \left( \lambda_n^{j'} \right)^2 t }{  \left( \lambda_n^j \right)^2} , \frac{ \lambda_n^{j'} x + x_n^{j'} - x_n^{j} }{ \lambda_n^j}  \right)  \right|^\frac4{d-2}   \left( \nabla \tilde{V}^{j'}  \right)(t,x)  \right\|_{L_{t,x}^\frac{2(d+2)}{d + 4} }  \xrightarrow{n\to\infty} 0.
		\end{align*}
		The proof of \eqref{eq4.52v32} is thus complete.
	\end{proof}
	
	Now, we are able to establish the existence and compactness up to scaling and translation of a critical element.
	
	\begin{theorem}[Existence and pre-compactness of a critical element] \label{pr5.2v32}
		There exists a maximal lifespan solution $u_c: [0, T^* ) \times \mathbb{R}^d \to \mathbb{C}$ to \eqref{eq1.1} with $T^* < \infty$  such that
		\begin{align*}
			E(u_c(0) ) = E_c , \ \|u_c(0) \|_{\dot{H}^1} < \|W\|_{\dot{H}^1},
			\text{ and }
			\|u_c \|_{S\left(  \left[ 0, T^*  \right) \right) } = \infty.
		\end{align*}
		Moreover, there are $x(t) \in \mathbb{R}^d $ and  $\lambda(t) \in \mathbb{R}_+$ for $t \in \big[0, T^*\big)$, such that
		\begin{align}\label{eq4.1v50}
			\mathcal{K} : =  \left\{ \frac1{\lambda(t)^{\f{d-2}2}} u_c  \left( t, \frac{ x - x(t) }{\lambda(t)}  \right): t \in [0, T^*)   \right\}
		\end{align}
		is pre-compact in $\dot{H}^1$.
	\end{theorem}
	
	\begin{proof}
		By Proposition \ref{pr4.6v32}, there exists at least one $j_0$ such that $\left\|v^{j_0}  \right\|_{ S(  [ 0, T^{j_0} ) )} = \infty.$
		Next, we notice that there is only one profile in the profile decomposition. Indeed, this follows from the definition of $E_c$, and the fact that the solutions of \eqref{eq1.1}, when the energy of the initial data is less than $E_c$ with kinetic energy of the initial data being less than the kinetic energy of $W$, has finite global Strichartz norm. In summary, we have
		\begin{align*}
			E \left(v^1(0)  \right) = E_c, \  \left\|v^1(0)  \right\|_{\dot{H}^1} < \|W\|_{\dot{H}^1},
			\left\|v^1 \right\|_{ S([0, T^1 ))  } = \infty,
			\text{ and }
			\left\|w_n^1  \right\|_{\dot{H}^1} \to 0, \ \text{as }n\to\infty.
		\end{align*}		
		Let $u_c : = v^1$ be a solution of \eqref{eq1.1} satisfying
		\begin{align*}
			E(u_c(0) ) = E_c \quad\text{ and }\quad \|u_c(0) \|_{\dot{H}^1} < \|W\|_{\dot{H}^1}.
		\end{align*}
		We shall show that $u_c(t)$ is pre-compact up to translation and scaling in $\dot{H}^1$ via a contradiction argument. If this fails, then there exist $\eta_0 > 0$ and a sequence $\{t_n\}_{n=1}^\infty$, $t_n \ge 0$ such that for any $\lambda_0 \in \mathbb{R}_+$, $x_0 \in \mathbb{R}^d$, we have
		\begin{align}\label{eq4.15v32}
			\left\|\frac1{\lambda_0^\frac{d-2}2} u_c \left(  t_n, \frac{x- x_0}{ \lambda_0} \right) - u_c ( t_m, x ) \right\|_{\dot{H}^1} \ge \eta_0, \ \text{for all}\, n \ne m.
		\end{align}
		Passing to a subsequence if necessary, we may assume that  $t_n \to \bar{t} \in [0, T^*]$.
		By continuity of the flow in $\dot{H}^1$, $\bar{t} = T^*$. Then $u_c(t+t_n)$ satisfies \eqref{eq5.2v32}.
		Hence, we can find a sequence $ \left(  x_n' , \lambda_n' \right) $ and {{}functions} $\phi \in \dot{H}^1$, $w_n \in \dot{H}^1$ such that
		\begin{align}\label{eq4.58v33}
			u_c(t_n, x ) = \frac1{  \left( \lambda_n'  \right)^\frac{d-2 }2 } \phi \left( \frac{ x- x_n'}{ \lambda_n'}  \right) + w_n ( x) ,
			\text{ with } \|w_n \|_{\dot{H}^1} \to 0, \text{ as } n \to \infty.
		\end{align}		
		By  \eqref{eq4.15v32} and \eqref{eq4.58v33}, for $ n \ne m$ large (independently of $\lambda_0$ and $x_0$), {{}we have}
		\begin{align*}
			\left\| \frac1{ \lambda_0^\frac{d-2}2} \frac1{ ( \lambda_n')^\frac{d-2}2} \phi \left(  \frac{ \frac{x-x_0}{ \lambda_0} - x_n'}{ \lambda_n'} \right)
			- \frac1{  \left( \lambda_m' \right)^\frac{d-2}2 } \phi \left( \frac{ x- x_m'}{ \lambda_m' } \right) \right\|_{\dot{H}^1} \ge \frac{\eta_0}2,
		\end{align*}
		and so
		\begin{align*}
			\left\|  \left( \frac{ \lambda_m'}{ \lambda_n' \lambda_0 } \right)^\frac{d-2}2  \phi \left( \frac{ \lambda_m' y  }{ \lambda_0 \lambda_n'} + \tilde{x}_{n,m} - \tilde{x}_0  \right) - \phi(y)  \right\|_{\dot{H}^1} \ge \frac{\eta_0}2,
		\end{align*}
		where $\tilde{x}_{n,m}$ is a suitable point in $\mathbb{R}^d$ and $\lambda_0, \tilde{x}_0$ are arbitrary.
		Taking $\lambda_0 = \frac{ \lambda_m'}{ \lambda_n'}$ and $\tilde{x}_0 = x_{n,m}$, we reach a contradiction. 		
	\end{proof}

	\section{Rigidity}\label{se6v23}
	
	In order to complete our proof of the global well-posedness for the case $\|u_0 \|_{\dot{H}^1} < \|W\|_{\dot{H}^1}$ in Theorem \ref{th1.1v3}, we need a rigidity theorem to exclude the critical element in Theorem \ref{pr5.2v32}.
	The key point is that dissipation stabilizes the behavior of solutions. The rigidity theorem for $d=3, 4$ is stated as follows.
	\begin{theorem}\label{thm:rigidity}
		The critical element $u_c$ in Theorem \ref{pr5.2v32} does not exist when $d = 3, 4$.	
	\end{theorem}
	
	The idea for the proof of rigidity is mainly motivated by \cite{KK11}, where the authors dealt with the incompressible Navier-Stokes equations, and it also relies on certain parabolic tools developed in \cite{ESS02, ESS031,ESS032}. The main difference between the proof of \cite{KK11} and ours is that in \cite{KK11}, they first constructed the linear profile decomposition for suitable weak solutions of the Navier-Stokes equations, and then applied the backward uniqueness and unique continuation for vorticity equations to exclude the minimal energy blowup solution; while in our case, for \eqref{eq1.1}, we can instead directly do this procedure on the equation itself.
	
	\subsection{Preliminary tools}
	
	There are three main ingredients to prove rigidity. The first one is the smoothness of solutions restricted on the parabolic domain $Q_r:= \left(-r^2, 0 \right) \times B_r(0)$ ($r>0$), with sufficiently small initial data. More precisely, we shall prove the following result.
	\begin{proposition}[Local smallness regularity criterion]\label{prop:regular}
		Set
		\begin{align}\label{est:theta}
			\epsilon_0 = ( \Re z)^\frac{d-2}8 .
		\end{align}
		{{}If} $u$ is a solution of \eqref{eq1.1} on $Q_1$ satisfying
		\begin{align}\label{eq5.2v50}
			\|u\|_{L^\infty_t \left(\dot{H}^1_x \cap L^{\f{2d}{d-2}}_x \right)(Q_1)}=\varepsilon < \varepsilon_0,
		\end{align}
		then {{}$u$ is smooth on $\overline{Q_{\frac12}}$ and for each $k \in \mathbb{Z}_+$, there exists a constant $C$ such that}
		\begin{align}\label{eq5.3vnew}
			\sup_{Q_{\frac12}} \left|\nabla^k u \right|\le C\varepsilon.
		\end{align}
	\end{proposition}
	\begin{proof}
		Define
		\begin{align*}
			\|u\|_{X(Q_1)}^2 := \|u\|_{L_t^\infty L_x^{\f{2d}{d-2}}(Q_1)}^2 + \|\nabla u\|_{L_t^2 L_x^{\f{2d}{d-2}}(Q_1)}^2 + \|\nabla u\|_{L_t^\infty L^2_x(Q_1)}^2 + \left\|\nabla^2 u \right\|_{L_{t,x}^2  (Q_1)}^2.
		\end{align*}
		Taking spatial derivatives on the first equation in \eqref{eq1.1}, we get
		\begin{align} \label{eq_diff1}
			\bar{z}(\nabla u)_t-\Delta (\nabla u) = \f{d}{d-2} |u|^{\f4{d-2}} \nabla u   +\f{2}{d-2}|u|^{\f{8-2d}{d-2}} (\nabla \bar{u})u^2.
		\end{align}
		Let $\phi_0$  be a smooth function such that {{}$0\le \phi_0 \le1$ and for some $\frac12< \rho_0 <1$},
		\begin{align*}
			\phi_0 (x) =
			\begin{cases}
				1, x\in B_{\rho_0}, \\
				0, x \not\in B_1,
			\end{cases}
		\end{align*}
		and $\phi_1$ be a smooth function such that for $\frac12<\rho_1<\rho_0$,
		\begin{align*}
			\phi_1(x) =
			\begin{cases}
				1, x \in B_{\rho_1}\\
				0, x \not\in B_{\rho_0}.
			\end{cases}
		\end{align*}
		Then $\phi_0\phi_1=\phi_1$ on the support of $\phi_1$.
		
		Taking the inner product with $\phi_0^2\nabla \bar{u}$ on both sides of \eqref{eq_diff1} and then integrating in space, we obtain
		\begin{equation*}
			\begin{aligned}
				\int_{|x|\le 1} \left[ \bar{z} (\nabla u)_t-\Delta (\nabla u) \right] \phi_0^2 \nabla \bar{u}  \, \mathrm{d}x
				=  \int_{|x|\le 1} \left(\f{d}{d-2}|u|^{\f4{d-2}} \nabla u  +\f{2}{d-2}|u|^{\f{8-2d}{d-2}}(\nabla \bar{u})u^2 \right) \phi_0^2  \nabla \bar{u}
				\,\mathrm{d}x .
			\end{aligned}
		\end{equation*}
		After taking real parts, we get
		\begin{equation*}
			\begin{aligned}
				& \int_{|x|\le 1} \left|\nabla^2u \right|^2\phi_0^2 \,\mathrm{d}x +  \frac{ \Re z}2  \frac{d}{dt}  \int_{|x|\le 1}\phi_0^2|\nabla u|^2 \,\mathrm{d}x \\
				&= 2 \int_{\rho_0 \le |x|\le 1} \Re \left[\nabla^2 u \nabla \phi_0 \nabla \bar{u} \phi_0 \right]\,\mathrm{d}x
				+\int_{|x|\le 1}\Re \left[\left(\f{d}{d-2}|u|^{\f4{d-2}}  \nabla u + \f{2}{d-2}|u|^{\f{8-2d}{d-2}} (\nabla \bar{u})u^2\right) \phi_0^2  \nabla \bar{u} \right]\,\mathrm{d}x.
			\end{aligned}
		\end{equation*}
		Integrating over $ \left[-\rho_0^2, t \right]$ with respect to $t \in  \left( - \rho^2, 0 \right)$, we then get
		\begin{equation}\label{local_2}
			\begin{aligned}
				& \frac{ \Re z }2 \int_{|x|\le 1}\phi_0^2|\nabla u (t)|^2 \,\mathrm{d}x + \int_{-\rho_0^2}^t\int_{|x|\le 1}\phi_0^2 \left|\nabla^2 u \right|^2 \,\mathrm{d}x \mathrm{d}s \\
				&  =\frac{ \Re z }2 \int_{|x|\le 1}\phi_0^2 \left|\nabla u \left(-\rho_0^2 \right) \right|^2 \,\mathrm{d}x
				+2\int_{-\rho_0^2}^t\int_{\rho_0 \le |x|\le 1} \Re \left[\nabla^2 u \nabla \phi_0 \nabla \bar{u} \phi_0 \right] \,\mathrm{d}x \mathrm{d}s \\
				&\quad +\int_{-\rho_0^2}^t\int_{|x|\le 1}\Re \left[\left(\f{d}{d-2} |u|^{\f4{d-2}}   \nabla u + \f{2}{d-2}|u|^{\f{8-2d}{d-2}} \nabla \bar{u} u^2\right) \phi_0^2 \nabla \bar{u} \right]\,\mathrm{d}x \mathrm{d}s .
			\end{aligned}
		\end{equation}
		By H\"older's inequality, we have
		\begin{equation}\label{local_3}
			\begin{aligned}
				\left|\int_{-\rho_0^2}^t\int_{\rho_0 \le |x|\le 1}\Re \left[\nabla^2 u \nabla \phi_0 \nabla \bar{u} \phi_0 \right] \,\mathrm{d}x  \mathrm{d}s \right|
				\lesssim
				\left\|\phi_0 \nabla^2 u \right\|_{L^2_{t,x}  }\|\nabla u\|_{L^2_{t,x} }
			\end{aligned}
		\end{equation}
		and
		\begin{align}\label{local_4}
			\left|\int_{-\rho_0^2}^t\int_{|x|\le 1}  \Re \left[\left(\f{d}{d-2}|u|^{\f4{d-2}} \nabla u  + \f{2}{d-2}|u|^{\f{8-2d}{d-2}}(\nabla \bar{u})u^2\right) \phi_0^2
			\nabla \bar{u}  \right] \,\mathrm{d}x \mathrm{d}s  \right|
			\lesssim \|u\|_{L^\infty_t L^{\f{2d}{d-2}}_x}^{\f4{d-2}}  \|\phi_0 \nabla u\|_{L^2_t L^{\f{2d}{d-2}}_x}^2.
		\end{align}
		Combining \eqref{local_2}, \eqref{local_3} and \eqref{local_4} all together, we deduce
		\begin{align*}
			\Re z  \|\phi_0\nabla u\|_{L^\infty_t L^2_x}^2 +  \left\|\phi_0 \nabla^2 u \right\|_{L^2_{t,x}   }^2
			\lesssim \|\phi_0\nabla u(-\rho_0^2)\|_{L^2_x}^2 + \|\nabla u\|_{L^2_{t,x}  }^2 + \|u\|_{L^\infty_t L^{\f{2d}{d-2}}_x}^{\f4{d-2}} \|\phi_0 \nabla u\|_{L^2_t L^{\f{2d}{d-2}}_x}^2,
		\end{align*}
		{{}which together with \eqref{eq5.2v50} and Sobolev's inequality implies}
		\begin{equation}\label{local_5}
			\begin{aligned}
				\Re z \|\phi_0\nabla u\|_{L^\infty_t L^2_x}^2 + \left\|\phi_0 \nabla^2 u \right\|_{L^2_{t,x}  }^2
				& \lesssim \varepsilon^2 + \|\nabla u\|_{L^\infty_tL^2_x}^2+ \|u\|_{L^\infty_t \dot{H}^1_x}^{\f4{d-2}} \|\phi_0 \nabla u\|_{L^2_t L^{\f{2d}{d-2}}_x}^2\\
				& \lesssim \varepsilon^2 \left(1+\|\phi_0 \nabla u\|_{L^2_t L^{\f{2d}{d-2}}_x}^2 \right).
			\end{aligned}
		\end{equation}
		For simplicity, we write $v_0=\phi_0\nabla u$. By Leibnitz's rule, we have
		\begin{align*}
			\|\nabla v_0\|_{L^2_{t, x}   }^2 \lesssim \|\nabla u\|_{L_t^\infty L_x^2}^2+ \left\|\phi_0 \nabla^2 u \right\|_{L^2_{t,x}   }^2,
		\end{align*}
		and then \eqref{local_5} along with Sobolev's inequality gives
		\begin{equation}\label{local_6}
			\begin{aligned}
				\Re z \|v_0\|_{L^\infty_t L^2_x}^2 + \|\nabla v_0\|_{L^2_{t, x}   }^2 +  \|v_0\|_{L^2_t L^{\f{2d}{d-2}}_x}^2
				\lesssim \varepsilon^2 \left(1+\|v_0\|_{L^2_t L^{\f{2d}{d-2}}_x}^2 \right).
			\end{aligned}
		\end{equation}
		Taking $\varepsilon_0$ sufficiently small, we arrive at the desired estimate
		\begin{align*}
			\| u\|_{X \left(Q_{\rho_0} \right)}\le \varepsilon.
		\end{align*}
		
		In order to apply the induction argument, we need to take higher order derivatives on \eqref{eq_diff1}.
		In the following, we consider separately the arguments for $d =3$ and $d = 4$.
		\medskip
		
		\textbf{Case I: $d=3$.}
		\medskip
		
		{{}Differentiating \eqref{eq_diff1} again}, we obtain
		\begin{equation}\label{eq_diff2}
			\begin{aligned}
				\bar{z}\left(\nabla^2 u \right)_t - \Delta \left(\nabla^2 u \right) & = 3 |u|^4  \nabla^2 u  + 6|u|^2 (\nabla u)(\nabla u) \bar{u}
				+ 6|u|^2  |\nabla u|^2 u \\
				&\quad+ 8|u|^2  |\nabla u|^2 u
				+4|u|^2 \left(\nabla^2 \bar{u} \right)u^2  +4 \left(\bar{u}\nabla u + u \nabla\bar{u} \right) (\nabla \bar{u})u^2.
			\end{aligned}
		\end{equation}
		Multiplying \eqref{eq_diff2} by $\phi_1^2 \nabla ^2\bar{u}$, integrating with respect to the spatial variables, and then taking the real parts gives
		\begin{equation}\label{local_7}
			\begin{aligned}
				& \  \frac{ \Re z }2 \frac{d}{dt}\int_{B_{\rho_0}} \phi_1^2 \left|\nabla^2 u \right|^2 \,\mathrm{d}x + \int_{B_{\rho_0}} \phi_1^2 \left|\nabla^3u \right|^2 \,\mathrm{d}x \\
				= & \int_{B_{\rho_0}} 2\nabla\phi_1 \phi_1 \left(\nabla^3 u \right)  \left(\nabla^2 \bar{u} \right) \,\mathrm{d}x
				+ \int_{B_{\rho_0}} 3|u|^4 \left(\nabla^2 u \right)\phi_1^2  \nabla^2 \bar{u}  \,\mathrm{d}x
				+ \int_{B_{\rho_0}} 6|u|^2\left[ \bar{u}(\nabla u)^2   +  u |\nabla u|^2 \right] \phi_1^2  \nabla^2 \bar{u}   \,\mathrm{d}x \\
				& + \int_{B_{\rho_0}} \left[ 8|u|^2 |\nabla u|^2  u
				+4|u|^2  \left(\nabla^2 \bar{u} \right)u^2\right] \phi_1^2  \nabla^2 \bar{u}   \,\mathrm{d}x + \int_{B_{\rho_0}}4 \left(\bar{u}\nabla u + u \nabla\bar{u} \right) \left(\nabla \bar{u} \right)u^2 \phi_1^2   \nabla^2 \bar{u}   \,\mathrm{d}x .
			\end{aligned}
		\end{equation}
		By \eqref{local_6}, there exists $t_1\in \left(-\rho_0^2, 0 \right)$ such that $\left\|\nabla v_0 (  t_1) \right\|_{ L^2 ( B_{\rho_0} ) } \lesssim \varepsilon$ (the implicit constant possibly depends on $\rho_0$). Set $v_1:=\phi_1 \nabla^2 u$. Then
		\begin{align*}
			\left\|v_1(t_1) \right\|_{L^2 (B_{\rho_0})} \le
			\|\nabla v_0 ( t_1)\|_{ L^2 ( B_{\rho_0} ) } +\|\nabla u\|_{L^\infty_t L^2_x} \lesssim \varepsilon.
		\end{align*}
		Integrating \eqref{local_7} with respect to time over $[t_1, t]$ for $t_1<t<0$, we obtain from the Cauchy-Schwarz and Sobolev inequalities that
		\begin{equation*}
			\begin{aligned}
				\Re z  \| v_1\|_{ L^\infty_t L^2_x}^2 + \left\| \phi_1 \nabla^3u \right\|_{L^2_{t,x }   }^2
				&\lesssim  \Re z \| v_1(t_1)\|_{ L^2_x}^2+ \left\| \phi_1 \nabla^3u \right\|_{L^2_{t, x}   } \left\|\nabla \phi_1 \nabla^2u \right\|_{L^2_{t, x}   }  + \|\phi_0 u\|_{L^\infty_t L^6_x}^4 \|v_1\|_{L^2_t L^6_x}^2\\
				&\quad + \|\phi_0 u\|_{L^\infty_t L^6_x}^3 \|\phi_1 v_0\|_{L^\infty_t L^6_x}
				\|v_0\|_{L^2_t L^6_x}\|v_1\|_{L^2_tL^6_x}\\
				&\lesssim \varepsilon^2+ \varepsilon^4\|v_1\|_{L^2_t L^6_x}^2+\varepsilon^4 \|\phi_1 v_0\|_{L^\infty_t L^6_x}\|v_1\|_{L^2_tL^6_x}
				+ \left\| \phi_1 \nabla^3u \right\|_{L^2_{t, x}   }  \left\|\nabla \phi_1 \nabla^2u \right\|_{L^2_{t, x}   }.
			\end{aligned}
		\end{equation*}
		Since
		\begin{align*}
			\nabla v_1=\phi_1\nabla^3 u + \nabla\phi_1 \nabla^2u,
		\end{align*}
		an application of the Cauchy-Schwarz inequality gives
		\begin{equation}\label{local_9}
			\begin{aligned}
				\Re z  \| v_1\|_{ L^\infty_t L^2_x}^2 + \| \nabla v_1\|_{L^2_{t, x}   }^2
				&\lesssim \varepsilon^2 + \varepsilon^4 \|v_1\|_{L^2_t L^6_x}^2+\varepsilon^4 \|\phi_1 v_0\|_{L^\infty_t L^6_x} \|v_1\|_{L^2_tL^6_x}
				+ \left\|\nabla \phi_1\nabla^2u \right\|_{L^2_{t, x}   }^2.
			\end{aligned}
		\end{equation}
		Thus, it remains to estimate the terms $\phi_1 v_0$ and $\nabla \phi_1\nabla^2u$.
		
		By the Sobolev inequality, the Leibnitz rule, and \eqref{local_6}, we have
		\begin{align}\label{local_10}
			\|\phi_1 v_0\|_{L^\infty_t L^6_x}\lesssim \|\phi_1 \nabla v_0\|_{L^\infty_t L^2_x}+\|\nabla\phi_1 v_0\|_{L^\infty_t L^2_x}\lesssim \|\phi_1 \nabla v_0\|_{L^\infty_t L^2_x}+\frac{\varepsilon}{\sqrt{ \Re z }}.
		\end{align}
		Since $v_0$ satisfies
		\begin{align}
			\nabla v_0 = \nabla \phi_0 \nabla u + \phi_0 \nabla^2 u,
		\end{align}
		it follows
		\begin{align}\label{local_11}
			\|\phi_1 \nabla v_0\|_{L^\infty_t L^2_x} \lesssim \|\phi_0\nabla u\|_{L^\infty_t L^2_x} + \|v_1\|_{L^\infty_t L^2_x}
			\lesssim \varepsilon + \|v_1\|_{L^\infty_t L^2_x}.
		\end{align}
		Note that in $B_{\rho_1}$, we have
		\begin{align}
			|\nabla \phi_1|  \left|\nabla^2u \right| \le\frac{|\nabla \phi_1|}{|\phi_0|}
			\left( |\nabla v_0| + |\nabla\phi_0 \nabla u| \right),
		\end{align}
		which further implies
		\begin{align}\label{local_12}
			\left\|\nabla \phi_1\nabla^2u \right\|_{L^2_{t, x} }^2 \lesssim \|\nabla v_0\|_{L^2_{t, x}   }^2 + \|\nabla u\|_{L^2_{t, x}   }^2
			\lesssim \|\nabla v_0\|_{L^2_{t, x} }^2 + \varepsilon^2.
		\end{align}
		Consequently, substituting \eqref{local_10}--\eqref{local_12} into \eqref{local_9}, we obtain
		\begin{equation*}
			\begin{aligned}
				\Re z \| v_1\|_{ L^\infty_t L^2_x}^2 + \| \nabla v_1\|_{L^2_{t , x}  }^2 &\lesssim \varepsilon^2 + \varepsilon^4 \|v_1\|_{L^2_t L^6_x}^2+\varepsilon^5 \left(1+\frac1{\sqrt{ \Re z  }} \right)\|v_1\|_{L^2_tL^6_x}\\
				&\qquad+\varepsilon^4 \| v_1\|_{ L^\infty_t L^2_x}\|v_1\|_{L^2_t L^6_x} + \|\nabla v_0\|_{L^2_{t, x}   }^2.
			\end{aligned}
		\end{equation*}
		Taking $\varepsilon_0$ small enough, {{}we conclude from} \eqref{est:theta}, the Cauchy-Schwarz inequality and \eqref{eq_diff2} that
		\begin{align}\label{eq_diff3}
			\Re z \| v_1\|_{ L^\infty_t L^2_x}^2 + \| \nabla v_1\|_{L^2_{t , x}  }^2+\|v_1\|_{L^2_t L^6_x}^2\lesssim \varepsilon^2 \left(1+ \varepsilon^8 +\frac{\varepsilon^8}{ \Re z } \right) \lesssim \varepsilon^2.
		\end{align}
		This yields
		\begin{align*}
			\left\|\nabla u \right\|_{X(Q_{\rho_1})}\lesssim \varepsilon.
		\end{align*}
		
		\textbf{Case II: $d=4$.}
		\medskip
		
		Differentiating \eqref{eq_diff1}, we have
		\begin{equation}\label{eq_diff2d4}
			\begin{aligned}
				\bar{z} \left(\nabla^2 u \right)_t - \Delta \left(\nabla^2 u \right) & = 2|u|^2  \nabla^2 u  + 2(\nabla u)(\nabla u) \bar{u}
				+ 2  |\nabla u|^2 u
				+ \left(\nabla^2 \bar{u} \right)u^2 + 2|\nabla u|^2  u.
			\end{aligned}
		\end{equation}
		Multiplying \eqref{eq_diff2d4} by $\phi_1^2 \nabla ^2\bar{u}$, integrating with respect to the spatial variables, and then taking the real parts gives
		\begin{equation}\label{local_7d4}
			\begin{aligned}
				& \ \frac{\Re z }2 \frac{d}{dt}\int_{B_{\rho_0}} \phi_1^2 \left|\nabla^2 u \right|^2 \,\mathrm{d}x + \int_{B_{\rho_0}} \phi_1^2|\nabla^3u|^2 \,\mathrm{d}x \\
				&= \int_{B_{\rho_0}} 2\nabla\phi_1 \phi_1 (\nabla^3 u) \left(\nabla^2 \bar{u} \right) \,\mathrm{d}x + \int_{B_{\rho_0}} 2|u|^2 \left(\nabla^2 u \right)\phi_1^2 \nabla^2 \bar{u}  \,\mathrm{d}x\\
				&  + \int_{B_{\rho_0}} \left[ 2\bar{u}(\nabla u)(\nabla u) + 2 u |\nabla u|^2\right] \phi_1^2   \nabla^2 \bar{u}  \,\mathrm{d}x
				+ \int_{B_{\rho_0}} \left[  \left(\nabla^2 \bar{u} \right)u^2 + 2 |\nabla u|^2 u\right] \phi_1^2   \nabla^2 \bar{u}   \,\mathrm{d}x .
			\end{aligned}
		\end{equation}
		As in the previous case, we may find $t_1\in  \left(-\rho_0^2, 0 \right)$ such that $ \left\|\nabla v_0 (   t_1) \right\|_{ L^2 \left( B_{\rho_0} \right) } \lesssim \varepsilon$ (the implicit constant possibly depends on $\rho_0$). Set $v_1:=\phi_1 \nabla^2 u$. Then
		\begin{align*}
			\left\|v_1(t_1)\right\|_{L^2 (B_{\rho_0})} \le \left\|\nabla v_0 (  t_1) \right\|_{ L^2 ( B_{\rho_0} ) }
			+\|\nabla u\|_{L^\infty_t L^2_x} \lesssim \varepsilon.
		\end{align*}
		Integrating \eqref{local_7d4} with respect to time over $[t_1, t]$ for $t_1<t<0$, we obtain from the Cauchy-Schwarz and Sobolev inequalities that
		\begin{equation*}
			\begin{aligned}
				& \  \Re z \| v_1\|_{ L^\infty_t L^2_x}^2 + \left\| \phi_1 \nabla^3u \right\|_{L^2_{t, x}  }^2\\
				&\lesssim \Re z \| v_1(t_1)\|_{ L^2_x}^2+ \left\| \phi_1 \nabla^3u \right\|_{L^2_{t , x}   }  \left\|\nabla \phi_1 \nabla^2u \right\|_{L^2_{t, x}   }  + \|\phi_0 u\|_{L^\infty_t L^4_x}^2 \|v_1\|_{L^2_t L^4_x}^2\\
				&\qquad+ \|\phi_0 u\|_{L^\infty_t L^4_x}  \|\phi_1 v_0\|_{L^\infty_t L^4_x} \|v_0\|_{L^2_t L^4_x}\|v_1\|_{L^2_tL^4_x}\\
				&\lesssim \varepsilon^2+\varepsilon^2\|v_1\|_{L^2_t L^4_x}^2+\varepsilon^2  \|\phi_1 v_0\|_{L^\infty_t L^4_x}\|v_1\|_{L^2_tL^4_x} + \left\| \phi_1 \nabla^3u \right\|_{L^2_{t, x}   }  \left\|\nabla \phi_1 \nabla^2u \right\|_{L^2_{t, x}  }.
			\end{aligned}
		\end{equation*}
		Since
		\begin{align*}
			\nabla v_1=\phi_1\nabla^3 u + \nabla\phi_1 \nabla^2u,
		\end{align*}
		it follows from the Cauchy-Schwarz inequality that
		\begin{equation}\label{local_9d4}
			\begin{aligned}
				\Re z  \| v_1\|_{ L^\infty_t L^2_x}^2 + \| \nabla v_1\|_{L^2_{t, x}   }^2
				&\lesssim \varepsilon^2 + \varepsilon^2 \|v_1\|_{L^2_t L^4_x}^2+\varepsilon^2 \|\phi_1 v_0\|_{L^\infty_t L^4_x}\|v_1\|_{L^2_tL^4_x}+
				\left\|\nabla \phi_1\nabla^2u \right\|_{L^2_{t, x}   }^2.
			\end{aligned}
		\end{equation}
		So it suffices to control the terms $\phi_1 v_0$ and $\nabla \phi_1\nabla^2u$.
		
		By the Sobolev inequality and \eqref{local_6}, we have
		\begin{align}\label{local_10d4}
			\|\phi_1 v_0\|_{L^\infty_t L^4_x}\lesssim \|\phi_1 \nabla v_0\|_{L^\infty_t L^2_x}+\|\nabla\phi_1 v_0\|_{L^\infty_t L^2_x}\lesssim \|\phi_1 \nabla v_0\|_{L^\infty_t L^2_x}+\frac{\varepsilon}{\sqrt{ \Re z }}.
		\end{align}
		Since $v_0$ satisfies
		\begin{align}
			\nabla v_0 = \nabla \phi_0 \nabla u + \phi_0 \nabla^2 u,
		\end{align}
		we have
		\begin{align}\label{local_11d4}
			\|\phi_1 \nabla v_0\|_{L^\infty_t L^2_x} \lesssim \|\phi_0\nabla u\|_{L^\infty_t L^2_x} + \|v_1\|_{L^\infty_t L^2_x}
			\lesssim \varepsilon + \|v_1\|_{L^\infty_t L^2_x} .
		\end{align}
		Note that in $B_{\rho_1}$, we have
		\begin{align}
			|\nabla \phi_1|  \left|\nabla^2u \right|  \le\frac{|\nabla \phi_1|}{|\phi_0|} \left( |\nabla v_0| + |\nabla\phi_0 \nabla u| \right),
		\end{align}
		which gives the bound
		\begin{align}\label{local_12d4}
			\left\|\nabla \phi_1\nabla^2u \right\|_{L^2_{t,x}  }^2 \lesssim \|\nabla v_0\|_{L^2_{t,x} }^2
			+ \|\nabla u\|_{L^2_{t,x}}^2 \lesssim \|\nabla v_0\|_{L^2_{t,x} }^2 + \varepsilon^2.
		\end{align}
		Consequently, we {{}can} substitute \eqref{local_10d4}--\eqref{local_12d4} into \eqref{local_9d4} and deduce
		\begin{equation*}\label{local_13d4}
			\begin{aligned}
				& \ \Re z \| v_1\|_{ L^\infty_t L^2_x}^2 + \| \nabla v_1\|_{L^2_{t, x}   }^2\\
				&\lesssim \varepsilon^2 + \varepsilon^2 \|v_1\|_{L^2_t L^4_x}^2+\varepsilon^3 \left(1+\frac1{\sqrt{ \Re z }} \right)\|v_1\|_{L^2_tL^4_x} +\varepsilon^2 \| v_1\|_{ L^\infty_t L^2_x}\|v_1\|_{L^2_tL^4_x} + \|\nabla v_0\|_{L^2_{t, x}   }^2.
			\end{aligned}
		\end{equation*}
		Taking $\varepsilon_0$ small enough, {{}we conclude from} \eqref{est:theta}, the Cauchy-Schwarz inequality and \eqref{eq_diff2} that
		\begin{align}\label{eq_diff3d4}
			\Re z \| v_1\|_{ L^\infty_t L^2_x}^2 + \| \nabla v_1\|_{L^2_{t,x }  }^2+\|v_1\|_{L^2_t L^4_x}^2
			\lesssim \varepsilon^2 \left(1+\frac{\varepsilon^4}{\Re z } \right)\lesssim \varepsilon^2.
		\end{align}
		This yields
		\begin{align*}
			\left\|\nabla u \right\|_{X(Q_{\rho_1})}\lesssim \varepsilon.
		\end{align*}
		
		Since the nonlinear terms are of lower order, we can iterate this procedure for finitely many times.
		Consequently, for any given $k\in \mathbb{Z}_+$, there exists $\varepsilon_0=\varepsilon_0(k)$ such that under \eqref{eq5.2v50}, it holds
		\begin{align*}
			\left\|\nabla^k u \right\|_{X \left(Q_{\frac12} \right)}\le C \varepsilon,
		\end{align*}
		where $C$ depends on $k$.
		This completes the proof.
	\end{proof}
	
	For the proof of Theorem \ref{thm:rigidity}, we shall need two more important results from the theory of parabolic operators: the backward uniqueness (Theorem \ref{thm:back}) and unique continuation (Theorem \ref{thm:continuation}); we refer the interested readers to \cite{ESS02, ESS031, ESS032} for more results on these. Note however that, in the original references, Theorems \ref{thm:back} and \ref{thm:continuation} were stated with constant $c_0=1$.
	Since the general result for $c_0 > 0 $ can be proved similarly with $\sqrt{c_0}\nabla$ in place of $\nabla$, we state them in the general form below.
	\begin{theorem}[Backward uniqueness]\label{thm:back}
		{{}For fixed $c_0$, $R$, $\delta >0$, set} $Q_{R, \delta} :=  (-\delta, 0) \times  \left( \mathbb{R}^d \backslash B(0, R) \right)$. If a function {$v\colon  Q_{R, \delta}\to \mathbb{R}$} satisfies
		\begin{itemize}
			\item[(i)] $v \in H^2(\Om)$ for any bounded $\Om \subseteq Q_{R, \delta}$; \\
			\item[(ii)] $|v(t, x )| \le e^{M|x|^2}$ for some $M > 0$ and for all $( t, x ) \in Q_{R, \delta}$; \\
			\item[(iii)] $ \left|\partial_t v - c_0  \Delta v \right| \le c_1  \left( |v| + |\nabla v| \right)$ on $Q_{R, \delta}$ for some $c_1 > 0$; \\
			\item[(iv)] $v(  0, x )=0$ for all $x\in \mathbb{R}^d \backslash B(0, R)$,
		\end{itemize}
		then $v \equiv 0$ in $Q_{R, \delta}$.
	\end{theorem}

	\begin{theorem}[Unique continuation]\label{thm:continuation}
		For fixed $c_0$, $R$, $\delta >0$, set $Q_{R, \delta} :=  (-\delta, 0) \times  \left( \mathbb{R}^d \backslash B(0, R) \right)$.
		If a function $v\colon  Q_{R, \delta}\to \mathbb{R}$ satisfies
		\begin{itemize}
			\item[(i)] $v \in H^2( Q_{R, \delta})$; \\
			\item[(ii)]  $\forall\, k \in \mathbb{N}$, there exist {$C_k>0$}, such that $ | v(t, x)|  \le C_k  \left( |x| + \sqrt{-t} \right)^k$ {for all $  x\in \mathbb{R}^d \backslash B(0, R)$}; \\
			\item[(iii)]{{}for some $c_1 > 0$}, $ \left|\partial_t v  - c_0 \Delta v \right| \le c_1 \left(|v| + |\nabla v| \right)$ a.e. on $Q_{R, \delta}$,
		\end{itemize}
		then $v(0, x ) \equiv 0$ in $B(0, R)$.
	\end{theorem}
	
	\subsection{Proof of Theorem \ref{thm:rigidity}}
	
	The strategy we shall use in the proof of Theorem  \ref{thm:rigidity} is as follows:
	we first show that the center $x(t)$ in \eqref{eq4.1v50} is bounded. Then we consider the parabolic domain away from $x(t)$ and use the local smallness regularity and zero $L^2$ integral limit of the critical element to {{}deduce that} the support of the critical element is bounded. Finally we use the backward uniqueness and unique continuation theorems to rule out the critical element. The zero $L^2$ integral limit of critical element shall be proved in Lemma \ref{lem:zero_limit}, using the existence of limit point of subsequence of profile decomposition from Lemma \ref{lem:convergenc} below. The boundedness of $x(t)$ {{}shall be proved} in Lemma \ref{lem:bound_x}, using the positivity of infimum energy of critical element from Lemma \ref{lem:energy_artifical}.
	
	For convenience, we first show the uniform boundedness of the sequence $ \{ v_n(x)\}$ in $\dot{H}^1$ for $t_n \in [0, T^\ast)$ with $t_n\nearrow T^*$ as $n\to \infty$, where
	\begin{align}\label{eq5.44newv45}
		v_n(x) : = \frac1{\lambda(t_n)^\frac{d -2}2 } u_c \left( t_n, \frac{x-x(t_n)}{\lambda(t_n)} \right).
	\end{align}
	For simplicity, we write $\lam_n = \lam (t_n)$ and $x_n = x(t_n)$.
	\begin{lemma}\label{lem:convergenc}
		There exists $v \in \dot{H}^1$ such that up to an extraction of subsequences,
		\begin{align}\label{eq5.28v45}
			v_n \xrightarrow{n\to\infty} v \ \text{  in }\ \dot{H}^1.
		\end{align}
	\end{lemma}
	
	\begin{proof}
		\textbf{Step 1}. We first show that $\liminf\limits_{t\nearrow T^\ast-} \lam(t) \sqrt{T^\ast-t} >0$.
		\medskip
		
		Note that this is nontrivial since $\lam (t) \to \infty$ as $t \to T^\ast$. Suppose for contradiction that there is a sequence $\{ t_n \}$ such that $t_n \nearrow T^\ast$. Then up to an extraction of subsequence, $v_n \to \hat {v}$ strongly in $\dot{H}^1$ by the compactness.
		
		Let $\bar{T}>0$ and $\bar{T}_n>0$ be the right endpoint of the maximal lifespan of solutions to \eqref{eq1.1} with initial data $\hat{v}(x)$ at $T^\ast$ and $v_n(x)$ at $t_n$, respectively. Then the continuous dependence on initial data implies $0< \bar{T} \le \liminf\limits_{n \to \infty} \left(\bar{T}_n - t_n \right)$. {{}On the other hand}, the scaling gives
		\begin{align*}
			\bar{T}_n - t_n = T^{max} (v_n) = \lam_n^2 T^{max} \left(u_c (t_n) \right) = \lam_n^2 \left(T^\ast - t_n \right) \xrightarrow{n\to\infty} 0,
		\end{align*}
		which is a contradiction.
		\medskip
		
		\textbf{Step 2}. Convergence.
		\medskip
		
		By Theorem \ref{pr5.2v32} 
		and the Sobolev embedding, for any $\varepsilon>0$, there exists $R_\varepsilon > 0 $ such that for any positive $t \le T^\ast$,
		\begin{align}\label{est:concentrate}
			\int_{|x - x(t)| \ge \f{R_\varepsilon}{\lam(t)}} |\nabla u_c (t, x )|^2 + |u_c (t, x )|^{\f{2d}{d-2}} \,\mathrm{d}x < \varepsilon,
		\end{align}
		which implies
		\begin{align*}
			\int_{|x|\ge R_\varepsilon} |\nabla v_n (t, x )|^2 \,\mathrm{d}x < \varepsilon.
		\end{align*}
		Consequently, there exists $v \in \dot{H}^1$ such that up to an extraction of subsequences,
		\[
		v_n \xrightarrow{n\to\infty} v \quad \text{  in }\ \dot{H}^1.
		\]
	\end{proof}
	
	Using Lemma~\ref{lem:convergenc}, we arrive at the following result.
	\begin{lemma} \label{lem:zero_limit}
		For any $R>0$,
		\begin{align}\label{eq5.29vnew}
			\lim_{n\to\infty} \int_{|x|\le R} |u_c(t_n, x )|^2 \,\mathrm{d}x=0.
		\end{align}
	\end{lemma}
	\begin{proof}
		By \eqref{eq5.44newv45}, we have
		\begin{equation}\label{eq_convergence_u1}
			\begin{aligned}
				\int_{|x|\le R} |u_c(t_n, x )|^2 \,\mathrm{d}x = \f1{\lam_n^2} \|v_n\|_{L^2( B(x_n, \lam_n R))}^2.
			\end{aligned}
		\end{equation}
		In order to show \eqref{eq5.29vnew}, 
		we fix $\varepsilon >0$ and split the $L^2$-norm of $v_n$ into {{}two parts}:
		\begin{align*}
			\|v_n\|_{L^2 \left( B \left(x_n, \lam_n R \right) \right)}^2 = \|v_n\|_{L^2 \left( B \left(x_n, \lam_n R \right) \cap B(0, \varepsilon \lam_n R) \right)}^2 + \|v_n\|_{L^2  \left( B \left(x_n, \lam_n R \right) \cap B \left(0, \varepsilon \lam_n R \right)^c \right)}^2.
		\end{align*}
		By H\"older's inequality and Lemma \ref{lem:convergenc}, we have
		\begin{equation} \label{est_convergenc_u1}
			\begin{aligned}
				\|v_n\|_{L^2( B(x_n, \lam_n R) \cap B(0, \varepsilon \lam_n R))}^2
				&\lesssim \varepsilon^2 \lam_n^2 R^2\left( \|v_n-v\|_{L^{\f{2d}{d-2}}(\mathbb{R}^d)}^2 + \|v\|_{L^{\f{2d}{d-2}}(\mathbb{R}^d)}^2 \right),
			\end{aligned}
		\end{equation}
		and
		\begin{equation}\label{est_convergenc_u2}
			\begin{aligned}
				\|v_n\|_{L^2( B(x_n, \lam_n R) \cap B(0, \varepsilon \lam_n R)^c)}^2
				& \lesssim \lam_n^2 R^2\left( \|v_n-v\|_{L^{\f{2d}{d-2}}(\mathbb{R}^d)}^2 + \|v\|_{L^{\f{2d}{d-2}}(B(0, \varepsilon \lam_n R)^c)}^2 \right).
			\end{aligned}
		\end{equation}
		Substituting \eqref{est_convergenc_u1} and \eqref{est_convergenc_u2} into \eqref{eq_convergence_u1}, we arrive at
		\begin{equation*}
			\begin{aligned}
				\int_{|x|\le R} |u_c(  t_n, x )|^2 \,\mathrm{d}x
				\lesssim R^2\left( \|v_n-v\|_{L^{\f{2d}{d-2}}(\mathbb{R}^d)}^2 + \|v\|_{L^{\f{2d}{d-2}}(B(0, \varepsilon \lam_n R)^c)}^2 \right) +\varepsilon^2 \|v\|_{L^{\f{2d}{d-2}}(\mathbb{R}^d)}^2.
			\end{aligned}
		\end{equation*}
		By \eqref{eq5.28v45} and the Sobolev inequality, we infer
		\[
		\lim_{n \to \infty} \|v_n-v\|_{L^{\f{2d}{d-2}}(\mathbb{R}^d)}^2 = 0.
		\]
		Letting $\lam_n \to \infty$, we conclude
		\[
		\lim_{n \to \infty}\|v\|_{L^{\f{2d}{d-2}} \left(B(0, \varepsilon \lam_n R)^c \right)}^2 = 0.
		\]
		Since $\varepsilon$ is arbitrary, \eqref{eq5.29vnew} 
		follows.
	\end{proof}
	
	We will show the center $x(t)$ of the critical element in Theorem \ref{pr5.2v32} is bounded. Before doing this, we first show that the energy of critical element has a strictly positive lower bound.
	\begin{lemma}\label{lem:energy_artifical}
		$ \underline{E}:= \inf\limits_{0 \le t < T^\ast} E(u_c(t)) >0 $.
	\end{lemma}
	\begin{proof}
		Let $\phi \in C_0^\infty \left(\mathbb{R}^d \right)$ be a radial decreasing cut-off function such that
		\begin{align*}
			\phi(x) =
			\begin{cases}
				1, |x| \le 1, \\
				0, |x | \ge 2,
			\end{cases}
		\end{align*}
		and $\phi_R(x) = \phi \left( \frac{x}R \right)$. {{}Then} for each $t\in [0, T^\ast)$,
		\begin{equation*}
			\begin{aligned}
				\f12 \f{d}{dt}\int |u_c(t, x )|^2 \phi_R(x)  \,\mathrm{d}x
				= \f12 \int \left[ \bar{z} \Delta u_c \overline{u_c} + z u_c \Delta \overline{u_c} \right] \phi_R \,\mathrm{d}x
				+ \int \Re z |u_c|^{\f{2d}{d-2}} \phi_R \,\mathrm{d}x .
			\end{aligned}
		\end{equation*}
		Integrating by parts gives
		\begin{align*}
			\f12 \f{d}{dt}\int |u_c(t, x )|^2 \phi_R(x) \,\mathrm{d} x
			= \int \Re z  \left( |u_c |^{\f{2d}{d-2}} - |\nabla u_c|^2 \right) \phi_R \,\mathrm{d}x
			- \f1R \int \Re z \ \Re (\nabla u_c \overline{u_c}) (\nabla \phi)_R \,\mathrm{d}x, 
		\end{align*}
		where $ ( \nabla \phi)_R(x) = ( \nabla \phi) \left( \frac{x}R \right)$.
		Sobolev's inequality implies that $ \left| \int \Re z \big( |u_c|^{\f{2d}{d-2}} - |\nabla u_c|^2 \big) \phi_R \,\mathrm{d}x \right|$ is uniformly bounded. By Hardy's inequality, we have
		\begin{align*}
			\left| \f1R \int \Re z \ \Re (\nabla u_c \overline{u_c}) (\nabla \phi)_R \,\mathrm{d}x \right| \lesssim \int \f{|u_c (x)|}{|x|}|\nabla u_c(x)| \,\mathrm{d}x \lesssim 1 .
		\end{align*}
		Consequently, there exists a constant $C>0$ such that
		\begin{align*}
			\left|\f12 \f{d}{dt}\int |u_c(t, x )|^2 \phi_R(x) \,\mathrm{d}x \right| \le C.
		\end{align*}
		Integrating with respect to time over $[t_0, t]$ for $t\le T^\ast$, we get
		\begin{align*}
			\left|\f12\int |u_c(t, x )|^2 \phi_R(x) \,\mathrm{d}x  - \f12\int |u_c( t_0, x )|^2 \phi_R(x) \,\mathrm{d}x \right| \le C (t - t_0).
		\end{align*}
		Passing to the limit $t \nearrow T^\ast$, we obtain from Lemma \ref{lem:zero_limit} that for any $R > 0$,
		\begin{align*}
			\left|\f12\int |u_c(  t_0, x )|^2 \phi_R(x) \,\mathrm{d}x \right| \le C (T^\ast - t_0),
		\end{align*}
		{{}from which we deduce $u_c(t_0) \in L^2$ by sending} $R \to \infty$. Since $t_0$ can be chosen freely,
		we further conclude that $u_c(t) \in L^2$ and
		\begin{align}\label{est:decrease}
			\f12 \int |u_c(t, x )|^2 \,\mathrm{d}x \le C (T^\ast - t).
		\end{align}
		Suppose now for contradiction that $\underline{E} > 0$ fails. {{}Then, we infer from \eqref{eq5.28v45} that}
		\begin{align*}
			\lim\limits_{n \to \infty} \left( 2E(v_n) -   \f12\int |v_n|^{\f{2d}{d-2}} \,\mathrm{d}x  \right)
			= 2 \underline{E} - \f12 \int |v|^{\f{2d}{d-2}} \,\mathrm{d}x  \le - \f12 \int |v|^{\f{2d}{d-2}} \,\mathrm{d}x  < 0.
		\end{align*}
		So by \eqref{eq5.11v15}, we would have
		\begin{align*}
			\lim_{n\to\infty}\f{d}{dt} \f12\int |v_n(t_n)|^2 \,\mathrm{d}x = \f{\Re z } 2 \int |v|^{\f{2d}{d-2}} \,\mathrm{d}x  >0,
		\end{align*}
		{{}which} implies that for all $ t\ \text{close to}\ T^\ast$,
		\begin{align*}
			\f{d}{dt} \f12\int |u_c(t, x )|^2 \,\mathrm{d}x > 0.
		\end{align*}
		Hence $\f12\int |u_c(t, x )|^2 \,\mathrm{d}x$ is increasing with respect to those $t$ that are close to $T^\ast$. But this clearly contradicts \eqref{est:decrease}. Thus it must hold $\underline{E} >0$.
	\end{proof}
	Based on Lemma \ref{lem:energy_artifical}, we now show the boundedness of $x(t)$ of the critical element in Theorem \ref{pr5.2v32}.
	\begin{lemma}\label{lem:bound_x}
		$\sup\limits_{0\le t < T^\ast} |x(t)| < \infty$.
	\end{lemma}
	\begin{proof}
		Suppose {{}for contradiction} that
		\begin{align}\label{eq5.66v45}
			|x(t_n)| \to +\infty \text{  as } t_n \nearrow T^\ast.
		\end{align}
		Let $\psi: [0, \infty) \to [0, 1]$ be a smooth cut-off function such that $\psi(r)=0$ when $r\le 1$ and $\psi(r)=1$ when $r\ge 2$, {{}and then} define $\psi_R(x):= \psi \left( \frac{ |x|}R \right)$ for $R>0$. For any $t_0 \in (0, T^\ast)$, there is $R_0 \ge 1$ such that
		\begin{align}\label{est:bound_x1}
			\int \left( \f12 |\nabla u_c(  t_0, x )|^2 - \f{d-2}{2d}|u_c(  t_0, x )|^{\f{2d}{d-2}} \right) \psi_{R_0}(x) \,\mathrm{d}x  \le \f14 \underline{E}.
		\end{align}
		Note that $\lam (t_n) \to +\infty$ as $n \to \infty$. By \eqref{eq5.66v45},
		for any $\varepsilon > 0$, there is $R_\varepsilon >0$ such that $ B\left( x(t_n), \frac{ R_{\varepsilon}}{\lam(t_n)} \right) \subseteq  B(0, 2R_0)^c$ when $n$ is sufficiently large. Thus, {{}we obtain from} \eqref{est:concentrate} that
		\begin{align*}
			\lim_{t\nearrow T^\ast}\int \left( \f12 |\nabla u_c(t, x )|^2 - \frac{d-2}{2d}|u_c(t, x )|^{\f{2d}{d-2}} \right) \psi_{R_0}(x) \,\mathrm{d}x  = \underline{E}.
		\end{align*}
		Hence, there exists $t_1 \in \left(t_0, T^\ast \right)$ such that
		\begin{align}\label{est:bound_x2}
			\int \left( \f12 |\nabla u_c (  t_1,x )|^2 - \f{d-2}{2d} |u_c(  t_1, x )|^{\f{2d}{d-2}} \right) \psi_{R_0}(x) \,\mathrm{d}x  \ge \f12 \underline{E}.
		\end{align}
		Combining \eqref{est:bound_x1} with \eqref{est:bound_x2} leads to
		\begin{align}
			\int \left[ \left( \f12 |\nabla u_c(  t_1, x )|^2 - \f{d-2}{2d} |u_c(  t_1, x )|^{\f{2d}{d-2}} \right)
			- \left( \f12 |\nabla u_c(  t_0, x )|^2 - \f{d-2}{2d}|u_c(  t_0, x )|^{\f{2d}{d-2}} \right) \right] \psi_{R_0}(x) \,\mathrm{d}x  \ge \f14 \underline{E}.
		\end{align}
		Note that $\left|\nabla\psi_{R_0} (x) \right| \lesssim \f1{R_0} \lesssim 1$ and so direct computation gives
		\begin{equation}\label{est:bound_x3}
			\begin{aligned}
				& \ \int \left[ \left( \f12 |\nabla u_c(  t_1, x )|^2 - \f{d-2}{2d} |u_c(  t_1, x )|^{\f{2d}{d-2}} \right) - \left( \f12 |\nabla u_c(  t_0, x )|^2 - \f{d-2}{2d}|u_c(  t_0, x )|^{\f{2d}{d-2}} \right) \right] \psi_{R_0}(x) \,\mathrm{d}x \\
				& = \int_{t_0}^{t_1} \f{d}{dt} \int\left( \f12 |\nabla u_c( t, x )|^2 - \f{d-2}{2d} |u_c(  t, x )|^{\f{2d}{d-2}} \right) \psi_{R_0}(x) \,\mathrm{d}x \\
				& = -\int_{t_0}^{t_1} \int \Re \left( \pa_tu_c \left(\overline{\Delta u_c}
				+ |u_c|^\frac4{d-2} \overline{u_c} \right) \right)\psi_{R_0} + \Re \left( \pa_t u_c \overline{\nabla u_c} \cdot \nabla\psi_{R_0}\right)  \,\mathrm{d}x \\
				& =-\int_{t_0}^{t_1}\int \Re z |\pa_tu_c|^2 \psi_{R_0}+ \Re \left( \pa_tu_c \overline{\nabla u_c} \cdot \nabla\psi_{R_0}\right) \,\mathrm{d}x
				\lesssim  \int_{t_0}^{t_1}\int | \pa_tu_c | |\nabla  {u_c}| \,\mathrm{d}x.
			\end{aligned}
		\end{equation}
		Applying H\"older's inequality together with the uniform boundedness of kinetic energy of $u_c$, we get
		\begin{equation}\label{est:bound_x4}
			\begin{aligned}
				\int_{t_0}^{t_1}\int | \pa_t u_c | |\nabla  {u_c}| \,\mathrm{d}x \mathrm{d}t
				\lesssim \|\pa_tu_c\|_{L_{t,x}^2([t_0, T^\ast) \times \mathbb{R}^d)}\|\nabla u_c\|_{L^\infty_t L^2_x} \sqrt{T^\ast-t_0} \lesssim \|\pa_tu_c\|_{L_{t,x}^2( [t_0, T^\ast) \times \mathbb{R}^d)}.
			\end{aligned}
		\end{equation}
		Finally, combining \eqref{est:bound_x2}--\eqref{est:bound_x4} and \eqref{eq2.3v3'}, we conclude
		\begin{align*}
			0 < \f14 \underline{E} \lesssim \|\pa_tu_c\|_{L^2([t_0, T^\ast) \times\mathbb{R}^d)} \to 0, \quad \text{as}\ t_0 \nearrow T^\ast.
		\end{align*}
		This is a contradiction and so the proof is complete.
	\end{proof}
	With all the above auxiliary results, we are now able to present the proof of Theorem \ref{thm:rigidity}.
	\begin{proof}[Proof of Theorem \ref{thm:rigidity}.]
		According to Lemma \ref{lem:bound_x}, $|x(t)|$ is bounded when $t \to T^\ast$. So we could choose a ball away from $x(t)$.
		By Theorem \ref{pr5.2v32}, for any given small $\varepsilon_0>0$, there exists $R>0$ large enough, such that for any $|x_0 |\geq R$, we have
		\begin{align*}
			\|u_c\|_{L^\infty_t \dot{H}^1_x \cap L^\infty_tL^{\f{2d}{d-2}}_x \left(\tilde{Q}_{T^\ast} \right)} < \varepsilon_0,
		\end{align*}
		where $\tilde{Q}_{T^\ast} =   \left(0, T^\ast \right) \times  B \left(x_0 , \sqrt{T^\ast} \right) $.
		
		By Proposition \ref{prop:regular}, $u_c$ is smooth in $\Om_{T^\ast} :=   \left[\frac34T^\ast, T^\ast \right] \times \left(\mathbb{R}^d\backslash B(0, R) \right) $ (up to appropriate scaling and shift) with $u_c$ being bounded as well as its derivatives of any order. In particular, $u_c(\cdot,x )$ is continuous at $T^\ast$ for any $x \in \mathbb{R}^d\backslash B(0, R)$.  Note also that by Lemma \ref{lem:zero_limit}, $u_c \left(  T^\ast, x \right) =0$ for any $x \in \mathbb{R}^d\backslash B(0, R)$. Hence, $v_c:=1- \cos (|u_c|^2)$ is smooth and bounded as well as its derivatives of any order, and $v_c\left(  T^\ast, x \right) =0=u_c(T^\ast,x)$ for any $x \in \mathbb{R}^d\backslash B(0, R)$. Replacing $u_c$ by $\frac{u_c}{M}$ for a large $M$ if necessary, we may assume that $|u_c|\leq \frac{1}{2}$ on $\Omega_{T^\ast}$.
		
		On the other hand, by a direct computation, we have
		\[
		\begin{aligned}
			(\pa_t - \Delta)v_c & = - \bar{u}_c (\pa_t u_c - \Delta u_c) \sin (|u_c|^2) - u_c (\pa_t \bar{u}_c - \Delta \bar{u}_c) \sin (|u_c|^2)\\
			& \quad + 2 |\nabla u_c|^2 \sin (|u_c|^2) + \cos (|u_c|^2)(u_c \nabla \bar{u}_c + \bar{u}_c \nabla u_c)^2\\
			& = - \bar{u}_c ((z-1)\Delta u_c + |u_c|^{2d/(d-2)} u_c) \sin (|u_c|^2) \\
			&\quad - u_c ( (\bar{z}- 1)\Delta \bar{u}_c + |u_c|^{2d/(d-2)} \bar{u}_c) \sin (|u_c|^2)\\
			& \quad + 2 |\nabla u_c|^2 \sin (|u_c|^2) + \cos (|u_c|^2)(u_c \nabla \bar{u}_c + \bar{u}_c \nabla u_c)^2.
		\end{aligned}
		\]
		Since $|u_c|\leq \frac{1}{2}$ on $\Omega_{T^\ast}$, we infer that
		\[
		\left|(\pa_t - \Delta)v_c \right| \le c_1 |u_c|^2 \le c_2 \Big|1- \cos (|u_c|^2) \Big|=c_2|v_c|.
		\]
		Applying Theorem \ref{thm:back} to $v_c$ gives $v_c \equiv 0$ and thus also $u_c \equiv 0$  in $\Om_{T^\ast}$.
		
		
		
		Set $\tilde{\Om}_{T^\ast} := \Big(\frac34 T^\ast, \frac78T^\ast \Big] \times \mathbb{R}^d  $. Then, $u_c$ and $v_c$ as well as their derivatives of any order belong to $L^2_{loc}(\tilde{\Om}_{T^\ast})$. {{}Since $v_c\equiv 0$ in $\Om_{T^\ast} \cap \tilde{\Om}_{T^\ast}$}, after an appropriate shift of local regions, applying Theorem \ref{thm:continuation}  to $v_c$, we get $v_c\equiv 0$ in $\tilde{\Om}_{T^\ast}$, which in return implies that $u_c \equiv 0$ in $\tilde{\Om}_{T^\ast}$. Finally, we can apply Theorem \ref{th2.1v3} to  {{}further conclude} $u_c \equiv 0$. This is a contradiction and so the proof of Theorem \ref{thm:rigidity} is complete.
	\end{proof}

	\section{Asymptotic decay of global solutions}\label{se4v23}
	
	In this section, we show the global solution $u$ of \eqref{eq1.1} must decay to zero in $\dot{H}^1$ as time tends to positive infinity when $E(u_0 ) < E(W)$ with $\|u_0  \|_{\dot{H}^1} < \|W\|_{\dot{H}^1}$. First of all, for the $H^1$ data, we shall prove the finiteness of $\|\nabla u \|_{L^2}$, by exploiting the $L^2$ dissipation relation. This in return allows us to reduce our problem to the small $\dot{H}^1$ data case and therefore deduce the finiteness of  $S(\mathbb{R}_+)$-norm of $u$.
	After this, to remove the extra assumption of $H^1$-boundedness of initial data,
	we split the initial data in frequency, and estimate a perturbed equation, which yields the finiteness of the $S(\mathbb{R}_+)$-norm.
	Finally, we show global solutions with finite $S(\mathbb{R}_+)$-norm must decay to zero in the $\dot{H}^1$ norm.
	
	\begin{theorem}\label{th3.1v3}
		If $u \in C_t^0 \dot{H}^1_x  \left( \mathbb{R}_+ \times  \mathbb{R}^d \right) $ is a solution to \eqref{eq1.1} satisfying
		\begin{align}\label{eq3.1v3}
			E(u_0 ) <E(W)\quad \text{and}\quad \|u_0  \|_{\dot{H}^1} < \|W\|_{\dot{H}^1},
		\end{align}
		then we have
		\begin{align}\label{eq6.2v51}
			\|u \|_{S(\mathbb{R}_+)} < \infty
			\intertext{ and }
			\lim\limits_{t \to \infty } \|u(t) \|_{\dot{H}^1 } = 0. \label{eq6.3v51}
		\end{align}
	\end{theorem}
	
	\begin{proof}
		
		\textbf {Step 1.} Finiteness of the $S(\mathbb{R}_+)$-norm.
		\medskip
		
		Suppose $u_0 \in H^1( \mathbb{R}^d )$. Then, {{}we claim that}  for some $\bar{\delta} > 0$,
		\begin{align}\label{eq6.4v74}
			\sup\limits_{t \ge 0 } \|u(t)\|_{L^2}^2 + 2 \bar{\delta} \|\nabla u \|_{L_{t,x}^2 (\mathbb{R}_+ \times \mathbb{R}^d )}^2 { \lesssim} \|u_0 \|_{L^2}^2.
		\end{align}
		In fact, by \eqref{eq3.1v3} and \eqref{eq2.11v3}, we have
		\begin{align*}
			\int |\nabla u(t)|^2 - |u(t)|^\frac{2d}{d-2}  \,\mathrm{d}x \ge \bar{\delta} \int |\nabla u(t)|^2 \,\mathrm{d}x
		\end{align*}
		and so
		\begin{align*}
			2 \bar{\delta} \|\nabla u \|_{L_{t,x}^2(\mathbb{R}_+ \times \mathbb{R}^d )}^2 \le 2 \int_0^\infty \int_{\mathbb{R}^d } |\nabla u(t)|^2 - |u(t)|^\frac{2d}{d-2}  \,\mathrm{d}x \mathrm{d}t
			\le \frac1{ \Re z} \|u_0 \|_{L_x^2}^2.
		\end{align*}
		{{}The claim \eqref{eq6.4v74} now follows from the above estimate together with the fact}
		\begin{align*}
			\sup\limits_{t \ge 0 } \|u(t) \|_{L^2} \le \|u_0 \|_{L^2} .
		\end{align*}
		By \eqref{eq6.4v74}, we have for any $\epsilon_0 > 0$, there exists some time $t_0 \in \mathbb{R}_+ $ such that
		$\|u(t_0) \|_{\dot{H}^1} \le \epsilon_0$. Thus we can directly apply Theorem \ref{th2.1v3}$(4)$
		to get \eqref{eq6.2v51}.
		
		To remove the extra assumption $u_0 \in L^2$, split $u_0 = w_0 + v_0$, with $\|w_0 \|_{\dot{H}^1} \ll 1$ and $v_0 \in H^1$. Let $w(t)$ be the solution of \begin{align*}
			\begin{cases}
				\bar{z}  \partial_t w = \Delta w +   f(w), \\
				w(0,x) = w_0(x).
			\end{cases}
		\end{align*}
		By Theorem \ref{th2.1v3}(4), 
		$w \in C_t^0 \dot{H}_x^1  \left( \mathbb{R}_+ \times \mathbb{R}^d \right)$ and
		\begin{align}\label{eq3.4v3}
			\|w \|_{L_{t,x}^\frac{2(d+2)}{d-2}   \left( \mathbb{R}_+ \times \mathbb{R}^d  \right)} + \|\nabla w \|_{L_t^\infty L_x^2 \cap L_{t,x}^\frac{2(d+2)}d  \left(\mathbb{R}_+ \times \mathbb{R}^d \right)} \lesssim \|\nabla w_0 \|_{L^2} \ll 1.
		\end{align}
		Set $v: = u - w$. Then it is a solution of
		\begin{align}\label{eq6.5newv45}
			\bar{z}  \partial_t v - \Delta v =   f(v+w) - f(w)
			= \mathcal{O} \left( v^\frac{d+2}{d-2} \right) + \mathcal{O}  \left( v^\frac4{d-2} w \right),
		\end{align}
		and when $d \ge 5$, the second term {{}disappears}.
		
		Multiplying \eqref{eq6.5newv45} by $\bar{v}$ and then integrating in space-time, we get
		\begin{align*}
			\|v(t) \|_{L^2}^2 - \|v_0 \|_{L^2}^2 + 2 \int_0^t \int |\nabla v|^2 \,\mathrm{d}x \mathrm{d}s
			= 2 \Re z  \Re \int_0^t \int \bar{v} \left( f(v+w) - f(w) \right) \,\mathrm{d}x \mathrm{d}s.
		\end{align*}
		By \eqref{eq3.4v3}, we have
		\begin{align*}
			\sup\limits_{t \ge 0 } \|\nabla v(t) \|_{L^2} < \|\nabla W \|_{L^2}.
		\end{align*}
		Hence  {{}it follows from \eqref{eq2.11v3} that} for some $\bar{\delta} > 0$,
		\begin{align*}
			\| v(t) \|_{L^2}^2 + \bar{\delta} \int_0^t \|\nabla v(s) \|_{L^2}^2 \,\mathrm{d}s \lesssim \|v_0 \|_{L^2}^2 + \int_0^t \int_{\mathbb{R}^d}  |v|
			\left( |v|^\frac4{d-2} |w| + |v| |w|^\frac4{d-2} \right)  \,\mathrm{d}x \mathrm{d}s.
		\end{align*}
		When $d \le 5 $, by H\"older's and Sobolev's inequalities, we have
		\begin{align*}
			\| v(t) \|_{L^2_x}^2 + \bar{\delta} \|\nabla v \|_{L_{t,x}^2}^2
			\lesssim &  \|v_0 \|_{L_x^2}^2  +  \left\| |v|^2 |w|^\frac4{d-2} \right\|_{L_{t,x}^1} +  \left\| |v|^\frac4{d-2} v w \right\|_{L_{t,x}^1}\\
			\lesssim  & \|v_0 \|_{L_x^2}^2  + \|w \|_{L_t^\infty L_x^\frac{2d}{d-2} }^\frac4{d-2} \|v\|_{L_t^2 L_x^\frac{2d}{d-2}}^2 + \| w \|_{L_t^\infty L_x^\frac{2d}{d-2} } \|v \|_{L_t^\infty L_x^\frac{2d}{d-2} }^\frac{6-d}{d-2} \|v\|_{L_t^2 L_x^\frac{2d}{d-2} }^2 \\
			\lesssim & \|v_0 \|_{L_x^2}^2 + \|\nabla w\|_{L_t^\infty L_x^2}^\frac4{d-2} \|\nabla v \|_{L_{t,x}^2}^2 + \|\nabla w \|_{L_t^\infty L_x^2} \|\nabla v \|_{L_t^\infty L_x^2}^\frac{6- d}{d-2} \|\nabla v \|_{L_{t,x}^2 }^2.
		\end{align*}
		When $d \ge 6$, the term $|v|^\frac4{d-2} v w$  is not included, and thus the last term on the right hand side disappears in higher dimensions.
		
		So by \eqref{eq3.4v3}, choosing $ \| \nabla w_0 \|_{L^2}$ small enough yields $\int_0^\infty \|\nabla v \|_{L_x^2}^2 \,\mathrm{d}t < \infty$. Hence there exists $T > 0$ for which $\| v (T) \|_{\dot{H}^1_x} < \| w_0 \|_{ \dot{H}_x^1}$, and so
		\begin{align*}
			\|\nabla u(T) \|_{L_x^2} \le 2 \|\nabla w_0 \|_{L_x^2}.
		\end{align*}
		Choosing $\|\nabla w_0 \|_{L^2}$ small enough,
		\eqref{eq6.2v51} follows from Theorem \ref{th2.1v3}(4).
		\medskip
		
		\textbf{Step 2}. Decay in $\dot{H}^1$ as $t\to \infty$.
		\medskip
		
		By the Strichartz estimate, continuity argument, and \eqref{eq6.2v51}, we have
		\begin{align}\label{eq6.3v50}
			\|\nabla u \|_{L_{t,x}^\frac{2(d+2)}d  ( \mathbb{R}_+ \times \mathbb{R}^d )} < \infty.
		\end{align}
		
		{{} Note that} the solution $u$ to \eqref{eq1.1} can be written as
		\begin{align*}
			u(t) = &  e^{tz \Delta }  u_0  +  \int_\tau^t e^{ (t - s  ) z \Delta }  z  f(u(s)) \,\mathrm{d}s
			+ \int_0^\tau e^{ (t - s) z \Delta }  z f(u(s))\,\mathrm{d}s
			: = I + II + III,
		\end{align*}
		for some $\tau$ to be determined later.
		
		For term $I$, approximating $\nabla u_0$ by $ v \in L^1 \cap L^2 $, and using \eqref{eq2.2v3}, we get
		\begin{align*}
			\|I \|_{\dot{H}^1}
			\le \| e^{t z \Delta } ( \nabla u_0 - v) \|_{L^2} + \| e^{ t z \Delta } v \|_{L^2}
			\le \|\nabla u_0 - v \|_{L^2} + \| e^{t z \Delta }  v \|_{L^2} \xrightarrow{t\to\infty} 0.
		\end{align*}
		We now treat term $II$, {{}where we shall fix} $\tau$. By \eqref{eq6.3v50}, for any $\epsilon > 0$, we can find $\tau$ such that
		\begin{align*}
			\|u \|_{L_{t,x}^\frac{2(d+2)}{d-2}  ([ \tau , \infty ) \times \mathbb{R}^d  )}
			+ \|\nabla u \|_{L_{t,x}^\frac{2(d+2)}d  ([\tau, \infty ) \times \mathbb{R}^d )} \le \epsilon.
		\end{align*}
		Since we are considering the limit {{}behaviour as} $t \to \infty$, we may assume $t > \tau \gg 1$. By the Strichartz estimate,
		\begin{align*}
			\| II\|_{\dot{H}^1} \lesssim \| u \|_{L_{t,x}^\frac{2(d+2)}{d-2} ( [ \tau,t ] \times \mathbb{R}^d )}^\frac4{d-2}  \|\nabla u \|_{L_{t,x}^\frac{2(d+2)}d  ( [\tau, t ] \times \mathbb{R}^d )} \lesssim \epsilon^\frac{d+2}{d-2} .
		\end{align*}
		Having fixed $\tau$ in this manner, we turn now to term $III$.
		Notice that
		\begin{align*}
			III = \int_0^\tau e^{ (t- s)z \Delta }  z  f(u(s))\,\mathrm{d}s
			= e^{ (t - \tau ) z \Delta }  \int_0^\tau e^{ (\tau - s) z \Delta }  z f(u(s))\,\mathrm{d}s .
		\end{align*}
		Since
		\begin{align*}
			\int_0^\tau e^{ (\tau - s) z \Delta }  z f(u(s))\,\mathrm{d}s \in \dot{H}^1,
		\end{align*}
		{{}using a} similar approximation argument as that for term $I$, we infer that
		\begin{align*}
			\| III \|_{\dot{H}^1} = \left\| e^{ (t - \tau ) z \Delta }  \int_0^\tau e^{ (\tau - s) z \Delta }  z f(u(s))\,\mathrm{d}s  \right\|_{\dot{H}^1} \xrightarrow{t\to\infty} 0.
		\end{align*}
		Since $\epsilon$ is arbitrary, \eqref{eq6.3v51} follows.
		
	\end{proof}

	\section{Proof of Theorem \ref{th1.2v21}}\label{se7v23}
	
	\subsection{Existence of global weak solutions}
	
	In this section, we show that the focusing energy-critical nonlinear Schr\"odinger equation has a global weak solution. To show this, we will study the inviscid limit of solutions of the Ginzburg-Landau equation and this relies heavily on the uniform boundedness of solutions to the Ginzburg-Landau equation.
	
	We first give an equivalent definition of weak-${H}^1$ solutions for nonlinear Schr\"odinger equation \eqref{eq1.1v18}; see for instance \cite{C} for the proof of equivalences.	
	\begin{definition}
		A complex-valued function $v$ on a time interval $I \subseteq \mathbb{R}$ is a weak-$H^1$ solution to \eqref{eq1.1v18} if $v\in L^\infty_t H_x^{1}(I \times \mathbb{R}^d)$, $\partial_tv \in L^\infty_t H_x^{-1}(I \times \mathbb{R}^d)$ such that \eqref{eq1.1v18} holds for almost every $t\in I$ in $H^{-1}_x$.
	\end{definition}
	
	The existence of global weak solutions to the defocusing energy-critical nonlinear Schr\"odinger equation was already proved in \cite{St}; see also \cite[Theorem 9.4.1]{C}. As was pointed out in Question \ref{ques:C} in the introduction of this paper, the existence of global weak solutions to the corresponding focusing energy-critical NLS in 3-D remains open for two decades. With the aid of Theorem \ref{th1.1v3}, we are now able to provide an affirmative answer to this open problem.
	
	\begin{proof}[Proof of Theorem \ref{th1.2v21} Part 1)]
		By Theorem \ref{th1.1v3}, there are global solutions $u_n$ of
		\begin{align}\label{eq7.1v24}
			\begin{cases}
				\partial_t u_n - z_n \Delta u_n - z_n f(u_n) = 0, \\
				u_n(0) = v_0,
			\end{cases}
		\end{align}
		where $|z_n| = 1$ with $\Re z_n \ge 0$ and $z_n\to i$ as $n \to \infty$.
		From Lemma \ref{le2.1v3}, we have $u_n \in L^\infty_t  H^1_x $ and $\pa_tu_n\in L^\infty_t H^{-1}_x$, with the uniform bound
		\begin{align}\label{est:uniform}
			\|u_n \|_{L^\infty_t H^1_x} + \left\|\pa_tu_n \right\|_{L^\infty_t H^{-1}_x} \le C \left(\|v_0\|_{H^1}\right).
		\end{align}
		Then there exists $v \in L^\infty_t {H}^1_x$ such that, up to an extraction of subsequence, we have
		\begin{align}\label{eq:weak_star}
			\begin{aligned}
				&(1)\  u_n(t) \rightharpoonup v(t)\ \text{weakly}\ \text{in}\ H^1_x ,\ \text{ as } n  \to \infty \ \text{for every} \ t \in \mathbb{R}_+ . \\
				& (2) \ \text{For every } t\in \mathbb{R}_+, \ 
				\  u_{n}(t,x)  \to v (t,x) \  \text{as} \ n \to \infty
				\text{ for almost every  } \ x\in \mathbb{R}^3,		\\
				& (3) \  u_n(t, x) \to v(t, x) \ \text{ as } n  \to \infty \ \text{for almost every} \ (t, x) \in \mathbb{R}_+ \times \mathbb{R}^3.
			\end{aligned}
		\end{align}
		Testing the equation \eqref{eq7.1v24} by $\varphi \in C_c^\infty(\mathbb{R}^3 )$ and $\psi \in  C^\infty_c( \mathbb{R}_+)$, we get
		\begin{align*}
			-\int_0^\infty \left( \int_{\mathbb{R}^3 } \bar{z}_n   u_n  \varphi  \,\mathrm{d}x \right) \psi_t\mathrm{d}t + \int_0^\infty \left(\int_{\mathbb{R}^3 } \nabla u_n  \cdot \nabla \varphi \,\mathrm{d}x\right) \psi \mathrm{d}t = \int_0^\infty \left(\int_{\mathbb{R}^3 } |u_n |^{4 } u_n  \varphi \,\mathrm{d}x\right) \psi\mathrm{d}t .
		\end{align*}
		Sending $n\to \infty$ and using the weak convergence in \eqref{eq:weak_star}, we infer that
		\begin{align*}
			\int_0^\infty \left(\int_{\mathbb{R}^3 } \bar{z}_n   u_n  \varphi \,\mathrm{d}x\right) \psi_t \mathrm{d}t
			\xrightarrow{n\to \infty} - \int_0^\infty \left(\int_{\mathbb{R}^3 } i v \varphi \,\mathrm{d}x \right) \psi_t\mathrm{d}t,
		\end{align*}
		and
		\begin{align*}
			\int_0^\infty \left(\int_{\mathbb{R}^3 } \nabla u_n  \cdot \nabla \varphi \,\mathrm{d}x \right) \psi\mathrm{d}t \xrightarrow{n\to \infty} \int_0^\infty \left(\int_{\mathbb{R}^3 } \nabla v \cdot \nabla \varphi \,\mathrm{d}x \right) \psi\mathrm{d}t.
		\end{align*}
		Note that $|u_n |^{4 } u_n  \varphi\psi$ is compactly supported and thus is bounded in $L^\infty_t L^{ \frac65}_x$. Then by the almost everywhere convergence in \eqref{eq:weak_star}, we get that $|u_n |^{4 } u_n  \varphi \psi \to |v|^{4 } v  \varphi \psi$ for almost every $(t, x)$, which  guarantees that $|u_n |^{4 } u_n  \varphi \psi \to |v|^{4 } v  \varphi \psi$ in $L^1_{t,x}$ as $n\to \infty$.
		
		Consequently, by taking $n \to \infty $, we obtain
		\begin{align*}
			\int_0^\infty\int_{\mathbb{R}^3 } i v \varphi \psi_t \,\mathrm{d}x \mathrm{d}t + \int_0^\infty\int_{\mathbb{R}^3 } \nabla v \cdot \nabla \varphi \psi\,\mathrm{d}x \mathrm{d}t  = \int_0^\infty\int_{\mathbb{R}^3 } |v|^{4 } v \varphi \psi\,\mathrm{d}x \mathrm{d}t .
		\end{align*}
		Hence $v$ satisfies
		\begin{equation}\label{eq7.2v87}
			iv_t +  \Delta v = - |v|^{4 } v
		\end{equation}
		in the sense of tempered distribution.
		
		By the lower semi-continuity of $L^2$-norm and the weak convergence in \eqref{eq:weak_star}, we have, for every $t\in\mathbb{R}_+$,
		\begin{align}\label{est:normdecrease_1}
			\|\nabla v(t)\|_{L^2}^2 \le \liminf_{n \to\infty}  \left\|\nabla u_{n } (t) \right\|_{L^2}^2.
		\end{align}
		Then Sobolev's inequality indicates that the sequence $u_{n}$ in \eqref{eq:weak_star}  satisfies $\int_{\mathbb{R}^3} \left|u_{n}(t,x) \right|^6 \le C$ for some positive constant $C$ uniformly for every $t\in\mathbb{R}_+$, so that the pointwise convergence of $u_{n}$ and Fatou's lemma give the estimate
		\begin{align}\label{est:normdecrease_2}
			\limsup_{n \to\infty} \int_{\mathbb{R}^3}  \left|u_{n} \right|^6 \mathrm{d}x \leq \int_{\mathbb{R}^3} |v|^6 \,\mathrm{d}x.
		\end{align}
		Combining \eqref{est:normdecrease_1} and \eqref{est:normdecrease_2}, we have
		\begin{align} \label{est:energy_nondecreas}
			\begin{aligned}
				E(v(t))
				\le \liminf_{n \to\infty} E \left(u_{n }(t) \right) \le E(v_0).
			\end{aligned}
		\end{align}
		For any fixed $T>0$, the map $t \mapsto \|v(t)\|_{L^2}$ is continuous on $[0,T]$, thus $v \in L^2_t \left([0, T] ; H^1_x \left(\mathbb{R}^3 \right) \right)$ and $\pa_tv \in L^2_t \left([0, T] ; H^{-1}_x \left(\mathbb{R}^3 \right) \right)$. By a standard limiting argument, we have
		\begin{align}
			\|v(t)\|_{L^2}^2 = \|v(s)\|_{L^2}^2 + \int_s^t \left< \pa_\tau v (\tau), v(\tau)\right>_{H^1, H^{-1}} \,\mathrm{d}\tau, \quad 0\le s \le t \le T,
		\end{align}
		which implies $v \in C \left([0, T]; L^2 \left(\mathbb{R}^3 \right) \right)$.
		
		Multiplying $i\bar{v}$ on both sides of \eqref{eq7.2v87} and integrating on $\mathbb{R}^3$, we obtain the mass conservation $M(v(t))= M(v_0)$ as desired.
		
		Thus, we have obtained a global weak solution of \eqref{eq1.1v18}
		on $[0, \infty)$. By time reflection, it is easy to extend the solution $v$ on $[0, \infty)$ to obtain a global weak solution $v$ on $\mathbb{R}$. 
	\end{proof}
	
	\subsection{Weak-strong uniqueness}
	In this subsection, we will show the weak-strong uniqueness when the $\dot{H}^1$-norm of the initial data is bounded by a constant which is related to the ground state $W$ and also the strong solution. First, we show that if \eqref{eq1.1v18} admits a smooth strong solution, then the weak solution is trapped around the strong solution in the following sense.
	\begin{lemma}
		\label{thm:stability}
		Suppose that $\tilde{v} \in  C^0_t H^2_x \left(\mathbb{R} \times \mathbb{R}^3 \right)$ is a strong solution to 3-D focusing NLS \eqref{eq1.1v18} with the initial data $\tilde{v}_0 \in H^2(\mathbb{R}^3)$ satisfying $E \left(\tilde{v}_0 \right) < E(W)$ and $ \left\|\nabla \tilde{v}_0 \right\|_{L^2} < \|\nabla W\|_{L^2}$. Let ${v}$ be the global weak solution to \eqref{eq1.1v18} with initial data $v_0$ satisfying the energy inequality
		\begin{align}
			E(v(t)) \le E(v_0), \quad \text{for all}\  t\in \mathbb{R}.
		\end{align}			
		Fix any $T>0$ and let $w = v-\tilde{v} $. Then there exists a constant $C= C  \left(\tilde{v}, W, T \right)>0$ such that
		\begin{align}
			\|w(t)\|_{H^1}^2 \le C  \left(1+\|w(0)\|_{H^1}^2 \right)
		\end{align}
		uniformly for all $ t \in (-T, T)$.
	\end{lemma}
	\begin{proof}
		Without loss of generality, we only consider the forward time case: $t\ge 0$. Direct computation shows that $w$ satisfies the equation
		\begin{align}\label{eq:perturbation}
			\begin{cases}
				iw_t + \Delta w = f \left(\tilde{v} \right) -   f \left(w + \tilde{v} \right) , \\
				w(0) = v_0- \tilde{v}_0 ,
			\end{cases}
		\end{align}
		where $f \left(\tilde{v} \right) = \left|\tilde{v} \right|^4 \tilde{v} $. We can expand $E(v)$ as
		\begin{align}\label{eq:decompose}
			E(v) = E \left(\tilde{v} \right) + I + II,
		\end{align}
		where
		\begin{align}
			I =  \Re \int_{\mathbb{R}^3} \left( \nabla \tilde{v} \cdot\overline{\nabla w }
			- f \left(\tilde{v} \right) \overline{w} \right) \,\mathrm{d}x ,
		\end{align}
		and
		\begin{align}
			II = \int_{\mathbb{R}^3} \left( \f12 |\nabla w|^2 - \left( F \left(w+ \tilde{v} \right) - F \left(\tilde{v} \right) - \Re \left( f \left(\tilde{v} \right) \bar{w}
			\right) \right) \right)
			\,\mathrm{d}x ,
		\end{align}
		with $F \left(\tilde{v} \right) = \f16 \left|\tilde{v} \right|^6$.
		
		Notice first that, by the Cauchy-Schwarz inequality,
		\begin{align}
			\f{d}{dt} \|w\|_{L^2}^2 =
			- 2  \Re \int_{\mathbb{R}^3}  \left( \overline{ f \left(\tilde{v} + w \right) - f \left(\tilde{v} \right) } \right) \cdot (iw)  \,\mathrm{d}x
			\le C_1 \int_{\mathbb{R}^3} \left(|w|^2 + F(w) \right) \,\mathrm{d}x .
		\end{align}
		By Gronwall's lemma and the fact that
		$\max
		\left\{
		\left\|\nabla \tilde{v} \right\|_{L^2}, \left\|\nabla v \right\|_{L^2} \right\} <  \|\nabla W\|_{L^2}$, we obtain
		\begin{align}\label{est:w_l2norm}
			\|w (t)\|_{L^2}^2 \le e^{C_1t} \|w (0)\|_{L^2}^2 + C_2 \int_0^t e^{C_1(t-s)} \|\nabla w(s)\|_{L^2}^2\,\mathrm{d}s,
		\end{align}
		where $C_1 = C_1 \left(\tilde{v} \right)$ and $C_2 = C_2(W)$.
		
		Then by a standard limitation argument and \eqref{est:w_l2norm}, we have
		\begin{align}\label{est:i_0-i_t}
			\begin{aligned}
				I(0) - I(t) &= - \int_0^t \f{d}{d s}  \left(\Re \int_{\mathbb{R}^3} \left( \nabla \tilde{v} \cdot\overline{\nabla w }
				- f \left(\tilde{v} \right) \overline{w} \right)(s,x)  \,\mathrm{d}x \right)
				\,\mathrm{d}s \\
				& \le C_3 \int_0^t \int_{\mathbb{R}^3}
				|w|^2 + F(w)  \,\mathrm{d}x \mathrm{d}s
				\le C_4 \|w (0)\|_{L^2}^2 + C_4 \int_0^t \|\nabla w(s)\|_{L^2}^2\,\mathrm{d}s
			\end{aligned}
		\end{align}
		with $C_4 = C_4 \left(T, \tilde{v}, W \right)$.
		
		Since $E \left(\tilde{v}(t) \right) = E \left(\tilde{v}_0 \right)$, and $E(v(t)) \le E(v_0)$, the expansion \eqref{eq:decompose} implies
		\begin{align}\label{est:v_t}
			0 \le E(v_0) - E(v(t)) = I(0) - I(t) + II(0) - II(t).
		\end{align}
		By the Cauchy-Schwarz inequality, Sobolev inequality and energy of ground state, we infer that
		\begin{align}\label{est:ii_t}
			\begin{aligned}
				II(t) &\ge \int_{\mathbb{R}^3} \left[ \f12|\nabla w|^2 - C F(w) - C_5 |w|^2  \right]  \,\mathrm{d}x
				\ge \f12\|\nabla w(t)\|_{L^2_x}^2  - C_6  - C_5 \|w(t)\|_{L^2}^2 ,
			\end{aligned}
		\end{align}
		where $C_5 = C_5 \left(\tilde{v} \right)$ and $C_6 = C_6 \left(\tilde{v}, W \right)$. At the same time, by the Cauchy-Schwarz inequality and Sobolev inequality, we have
		\begin{align}\label{est:ii_0}
			II(0) \le C_7  \left(\|\nabla w(0)\|_{L^2}^2 + \|w(0)\|_{L^2}^2 \right)
		\end{align}
		with $C_7 = C_7 \left(\tilde{v}, W \right)$
		
		We note that \eqref{est:ii_t} gives
		\begin{align}
			\|\nabla w(t)\|_{L^2}^2\le C_8  \left(II(t) + \|w(t)\|_{L^2}^2 + C(W)\right).
		\end{align}
		Then by \eqref{est:v_t}, \eqref{est:i_0-i_t} and \eqref{est:ii_0}, we finally conclude that
		\begin{align}\label{est:gronwall}
			\begin{aligned}
				\|\nabla w(t)\|_{L^2}^2  + \|w(t)\|_{L^2}^2 & \le C_9  \left(I(0) - I(t) + II(0) + \|w(t)\|_{L^2}^2  \right)\\
				& \le C \left(\|\nabla w(0)\|_{L^2}^2 + \|w(0)\|_{L^2}^2 +1 \right) + \int_0^t \|\nabla w(s)\|_{L^2}^2 + \|w(s)\|_{L^2}^2 \,\mathrm{d}s
			\end{aligned}
		\end{align}
		with a uniform constant $C= C \left(\tilde{v}, T, W \right)$.
	\end{proof}
	
	The above trapped result clearly implies the following short time stability result.
	\begin{corollary}
		Suppose that $\tilde{v} $ is the global strong solution to
		\eqref{eq1.1v18} with the initial data $v_0 \in H^2(\mathbb{R}^3)$ satisfying $E(v_0) < E(W)$ and $\|\nabla v_0\|_{L^2} < \|\nabla W\|_{L^2}$.
		Let ${v}$ be the global weak solution to \eqref{eq1.1v18} with initial datum $v_0$, satisfying the energy inequality \eqref{eq1.4v26}.
		Then for every $\varepsilon>0$, there exists a positive $T(\varepsilon)$, such that for each $t\in[0, T(\varepsilon))$, there holds
		\begin{equation}\label{eq:stability}
			\left\|v(t) - \tilde{v}(t) \right\|_{L^2}^2 < \varepsilon.
		\end{equation}
	\end{corollary}
	\begin{proof}
		From the trapped result in Lemma \ref{thm:stability} and \eqref{est:w_l2norm}, we obtain that for $w=v-\tilde{v}$,
		\begin{align} \label{est:w_2norm}
			\|w (t)\|_{L^2}^2 \le e^{C_1t} \|w (0)\|_{L^2}^2 + C \int_0^t e^{C_1(t-s)}  \left(1+\|w(0)\|_{H^1}^2 \right)\,\mathrm{d}s,
		\end{align}
		for any $0<  t \le T <\infty$. For $\|w(0)\|_{H^1}^2 = 0$, \eqref{est:w_2norm} reduces to
		\begin{align*}
			\|w (t)\|_{L^2}^2 \le  C(e^{C_1 t} -1).
		\end{align*}
		For each $\varepsilon>0$, we may choose $T = \frac{1}{C_1} \log \left(1 + C^{-1} \varepsilon \right)$ so that \eqref{eq:stability} holds for any $t \in (0, T)$.
	\end{proof}
	Now, we are able to prove the weak-strong uniqueness result of Theorem \ref{th1.2v21}.
	
	\begin{proof}[Proof of Theorem \ref{th1.2v21} Part 2)]
		Throughout the proof, the constant $C>0$ is allowed to change from line to line. We follow the proof of Lemma \ref{thm:stability} from \eqref{eq:perturbation} to \eqref{est:gronwall}, except that \eqref{est:ii_t} is replaced by
		\begin{align*}
			II(t) &\ge \int_{\mathbb{R}^3} \left[ \f12|\nabla w|^2 - C F(w) - C_5 |w|^2  \right]  \,\mathrm{d}x
			\ge \f12\|\nabla w(t)\|_{L^2_x}^2  - C_6\|w(t)\|_{L^6}^6  - C_5 \|w(t)\|_{L^2}^2 .
		\end{align*}
		Then we have 		\begin{align}\label{est:continuity}
			\begin{aligned}
				\|\nabla w(t)\|_{L^2}^2  + \|w(t)\|_{L^2}^2 \le C \left(\|w(t)\|_{L^6}^6 + \int_0^t \|\nabla w(s)\|_{L^2}^2 + \|w(s)\|_{L^2}^2 \,\mathrm{d}s \right),
			\end{aligned}
		\end{align}
		where we used the fact that $w(0) = 0$.
		
		We now turn to estimate $\|w(t)\|_{L^6}^6$ in \eqref{est:continuity}.
		Since $v$ is a weak solution to NLS, we may use the Sobolev embedding inequality,  \eqref{est:normdecrease_1}, variational results in Lemma \ref{le2.1v3} for solutions to Ginzburg-Landau equation, and second inequality in \eqref{est:energy_nondecreas}, to deduce that
		\begin{align}\label{est:w_2}
			\|v(t)\|_{L^6} \le C \left\|\nabla v_0 \right\|_{L^2}.
		\end{align}
		As $\tilde{v}$ is a strong solution to NLS, we employ the energy conservation, variational results of the solution to NLS and the Sobolev embedding inequality to get
		\begin{align}\label{est:w_3}
			\left\|\tilde{v}(t) \right\|_{L^6} \le C \|\nabla v_0\|_{L^2}.
		\end{align}
		Hence, by \eqref{est:w_2}, \eqref{est:w_3} and the triangle inequality, we have
		\begin{align*}
			\|w(t)\|_{L^6}^6  \le  C \left(\|v(t)\|_{L^6} +  \left\|\tilde{v}(t) \right\|_{L^6} \right)^4 \|\nabla w(t)\|_{L^2}^2  \le  C \left\|\nabla v_0 \right\|_{L^2}^4 \|\nabla w(t)\|_{L^2}^2.
		\end{align*}		
		From the condition of initial data $v_0$, \eqref{est:continuity} is reduced to
		\begin{align*}
			\begin{aligned}
				\|\nabla w(t)\|_{L^2}^2  + \|w(t)\|_{L^2}^2 \le \frac{C}{1-C\|\nabla v_0\|_{L^2}^4} \int_0^t \|\nabla w(s)\|_{L^2}^2 + \|w(s)\|_{L^2}^2 \,\mathrm{d}s,
			\end{aligned}
		\end{align*}
		which gives $w\equiv 0$ by Gronwall's inequality.
		
	\end{proof}

	\section{Remarks on the Galerkin method}\label{sec:concluding remark}

	The Gelerkin method is a powerful technique in proving the well-posedness theory for many (linear) partial differential equations; see for instance \cite[Chapter 2.5]{Li-Chen-1992-Book}. In this section, we shall point out that one cannot prove Theorem \ref{th1.2v21} easily via the classical Galerkin method. In the outline below, we shall discuss the major mathematical difficulties arsing from the Galerkin method.
	\medskip
	
	\textbf{Step 1.} Find a sequence of increasing subspace $V_N\subset V_{N+1}\subset \dot{H}^1$ so that their union is dense.
	
	In this step, we can make use of the Hermitian functions, that is, an orthonormal basis $\{h_n\}$ of $L^2(\mathbb{R})$ satisfying
	\begin{align}\label{hermitian}
		h_n ^\prime (x) + x h_n(x) = 2n h_{n-1}(x), \quad x \in  \mathbb{R},\\
		h_n ^\prime (x) - x h_n(x) = - h_{n+1}(x), \quad x \in  \mathbb{R},\notag\\
		h_n '' (x) - x^2 h_n(x) + (2n+1) h_n(x) =0, \quad x \in  \mathbb{R}.\notag
	\end{align}
	Then we can proceed to construct an orthonormal basis $\{w_j\}_{j=1}^\infty$ of $L^2(\mathbb{R}^3)$ through Hermitian functions.
	
	In order to use the Galerkin method, we need to find an orthonormal basis of $H^1(\mathbb{R}^3)$. From \eqref{hermitian}, $\{h_n ^\prime\}$ are not necessarily mutually orthogonal. Thus {some extra little work} is needed in order to find an orthonormal basis of $H^1(\mathbb{R}^3)$, which is defined on the entire (unbounded) space $\mathbb{R}^3$.
	\medskip
	
	\textbf{Step 2.} Construct local approximate solutions.

	Suppose that we were able to find an orthonormal basis $\{w_j\}_{j=1}^\infty$ of $H^1(\mathbb{R}^3)$. Then we define the finite-dimensional spaces
	\[
	V_N = \text{span} \{ w_1, \ldots, w_N\}.
	\]
	Given $T>0$, we want to find $u^N(t) = \sum_{j=1}^N d^N_j(t) w_j$, $t\in [0, T]$, solving
	\[
	i\int_{\mathbb{R}^3}u^N_t w +  \int_{\mathbb{R}^3}\nabla u^N \cdot \nabla w= \int_{\mathbb{R}^3}|u^N|^4u^N w
	\]
	for each $w \in V_N$. This amounts to solve the ODE
	\[
	i \frac{ \mathrm{d} \left(  {d}^N_j (t) \right) }{ \mathrm{d} t}(w_j, w_k)_{L^2} + d^N_j(t) (\nabla w_j, \nabla w_k)_{L^2} = (d^N_j(t))^5 (|w_j|^4 w_j, w_k)_{L^2}.
	\]
	To simplify our notation, we write $u_0^N=P_Nu_0$.
	
	If $u_0 \in \dot{H}^1$ satisfies $\|\nabla u_0 \|_{L^2} < \|\nabla W \|_{L^2}$ and $E(u_0 ) < (1 - \delta) E(W)$, then $u_0^N$ satisfies  $\|\nabla u_0^N \|_{L^2} < \|\nabla W \|_{L^2}$ and $E(u_0^N) < (1 -  \frac12 \delta) E(W)$. By the variational estimate as in Kenig and Merle \cite[Lemma 3.4]{KM1},  there exists $\bar{\delta} = \bar{\delta}  ( \delta, d ) > 0$ such that for all $t \in [-T,T]\subset (-T_{\min}, T_{\max} )$, the solution $u^N $ of focusing NLS 
	satisfies
	\begin{align*}
		& \int |\nabla u^N (t)|^2 \,\mathrm{d}x \le  \left( 1- \bar{\delta } \right) \int |\nabla W|^2\, \mathrm{d}x.
	\end{align*}
	
	\medskip

	\textbf{Step 3.} The Kenig-Merle arguments in dimension three.
	
	Now, the key technical step is to show that $T_{\min}=T_{\max}=\infty$. In dimension $d=3$, Kenig and Merle proved this {in the radial case}. Namely, they proved that the local strong solution $u \in C_t^0 ( (- T_{\min}, T_{\max}), \dot{H}_{rad}^1( \mathbb{R}^3)) $ can be extended to the global solution $u \in C_t^0 ( (- \infty, \infty), \dot{H}_{rad}^1( \mathbb{R}^3))$.
	
	The proof uses a contradiction argument and contains the following two main steps:
	\begin{enumerate}
		\item By the concentration-compactness argument \cite[Proposition 4.1]{KM1}, there exists a critical element $u_c$ so that the set
		\begin{equation}\label{eq:def for K}
			K =  \left\{ v: v(x, t) = \frac1{ \lambda(t)^\frac{1}2} u_c\left( \frac{x - x(t)}{\lambda(t)}, t\right) \right\}
		\end{equation}
		is compact in $\dot{H}^1$;
		
		\item By the rigidity theorem \cite[Theorem 5.1]{KM1}, such a critical element $u_c$ does not exist.
	\end{enumerate}
	For this argument to work, it is important that the $x(t)$ appearing in \eqref{eq:def for K} has to satisfy the decay estimate $|x(t) | = o(t)$ as $|t| \to \infty$. In dimension $d=3$, the radial condition ensures that $x(t)=0$ and thus in particular the above decay estimate is satisfied. {In the non-radial case when $d=3$, there is no known method to derive such a decay estimate} and it remains as a difficult open problem (see however the works by Killip-Visan \cite{KV2010} and Dodson \cite{D2019} in dimensions $d\geq 4$). Thus it is not easy to prove the global existence in the non-radial case.
	
	\medskip

	\textbf{Step 4.} Standard limiting argument.
	
	Applying a standing limit argument (up to a further subsequence), we find that for any given positive number $T$, $u^N \to u \in L^\infty([0, T]; H^1(\mathbb{R}^3))$. This however {does not mean that $u$ is a global weak solution} according to our definition. Note that as $u$ is not necessarily continuous in $t$ and we do not necessarily have uniqueness result for weak solutions of NLS, it is a {nontrivial task} to show that $u$ is indeed a global weak solution (according to our definition). In fact, even in the relatively simpler case of focusing mass-critical NLS, the proof of this is {far from trivial}; see the deep work of Terence Tao \cite{Tao-2009-APDE}.
	\bigskip

	\noindent \textbf{Acknowledgments. } 
	X. Cheng has been partially supported by the NSF of Jiangsu Province (Grant No.~BK20221497). C.-Y. Guo is supported by the Young Scientist Program of the Ministry of Science and Technology of China (No.~2021YFA1002200), the National Natural Science Foundation of China (No.~12101362 and 12311530037), the Taishan Scholar Project, the Natural Science Foundation of Shandong Province (No.~ZR2022YQ01) and the Jiangsu Provincial Scientific Research Center of Applied Mathematics (Grant No. BK20233002). Y. Zheng was supported by National Natural Science Foundation of China under Grant No.~11901350 and No.~12371172.
	
	The authors would like to thank Prof.~Gongbao Li for his constant encouragement in carrying out this project and this work is dedicated to celebrate his 70th birthday. They are grateful to Qingtang Su, Zaiyun Zhang and Lifeng Zhao for helpful discussion and to Terence Tao for kindly explaining some ideas of his interesting paper.

\end{document}